\documentclass[11pt]{amsart}
\usepackage[OT2,T1]{fontenc}
\usepackage[T1]{fontenc}
\usepackage[latin1]{inputenc}
\usepackage{a4}
\usepackage{setspace}
\usepackage{xypic}
\usepackage{amssymb}
\usepackage{amscd}
\usepackage{amsthm}
\usepackage[english]{babel}
\usepackage{latexsym}
\date{}
\newtheorem{thm}{Theorem}[section]
\newtheorem{lem}[thm]{Lemma}

\newtheorem{defn}[thm]{Definition}
\newtheorem*{thm*}{Theorem}
\newtheorem{rem}[thm]{Remark}

\newtheorem{Hypothesis}[thm]{Hypothesis}
\newtheorem{prop}[thm]{Proposition}

\newtheorem{cor}[thm]{Corollary}
\newtheorem{notation}[thm]{Notation}
\newenvironment{f-proof}[1][\sc D\'emonstration.]{\begin{trivlist}
\item[\hskip \labelsep {\bfseries #1}]}{\hfill{$\square$}\end{trivlist}}

\newcommand{\fonc}[5]{
 \begin{array}{cccc}
 #1: & #2 & \longrightarrow & #3\\
     & #4 & \longmapsto & #5
 \end{array}
}
\newcommand{\appl}[4]{
 \begin{array}{cccc}
  #1 & \longrightarrow & #2\\
  #3 & \longmapsto & #4
 \end{array}
}

\newcommand{\bq}{\mathbb Q}

\newcommand{\Z}{\mathbb Z}

\newcommand{\bL}{\mathbb L}
\newcommand{\br}{\mathbb R}

\newcommand{\X}{{\mathcal X}}

\newcommand{\et}{{\mathrm{et}}}
\newcommand{\red}{{\mathrm{red}}}
\newcommand{\Pic}{{\operatorname{Pic}}}

\title[Pontryagin duality for varieties over $p$-adic fields]
{Pontryagin duality for varieties over $p$-adic fields}

\date{}
\author{Thomas H. Geisser}
\address{
	Department of Mathematics, Rikkyo University, Ikebukuro, Tokyo, Japan
}
\email{geisser@rikkyo.ac.jp}
\author{Baptiste Morin}
\address{
	Department of Mathematics, Universit\'e de Bordeaux,
	Bordeaux, France
}
\email{Baptiste.Morin@math.u-bordeaux.fr}
\thanks{The first named author is Supported by JSPS Grant-in-Aid (C) 18K03258,
and the second named author by grant ANR-15-CE40-0002}

\begin{document}


\begin{abstract}
We define cohomological complexes of locally compact abelian groups associated with varieties over $p$-adic fields and prove a duality theorem under some assumption. Our duality takes the form of Pontryagin duality between locally compact motivic cohomology groups.
\end{abstract}


\maketitle


\section{Introduction}

Let $K$ be a finite extension of $\mathbb{Q}_p$. Let $\mathcal{O}_K$ be the ring of integers in $K$ and let $\X$ be a  regular, proper and flat scheme over $\mathcal{O}_K$
of dimension $d$. 
We denote by $\X_K$ its generic fiber and by $i:\X_s\rightarrow\X$ its special fiber. 
It is a classical result that for any integer $m>0$ we have a perfect duality of
motivic cohomology with finite coefficients 
$$H^{i}_{\et}(\X_K,\Z/m(n))\times 
H^{2d-i}_{\et}(\X_K,\Z/m(d-n))\rightarrow H^{2d}_{\et}(\X_K,\Z/m(d))\rightarrow \Z/m.$$
However, this does not lift to a duality of integral groups 
\cite{Geisser16}. For example, even for a curve $\X_K$, the dual of 
$H^1_{\et}(\X_K,\bq/\Z)\cong H^2_{\et}(\X_K,\Z)$ has both contributions from 
$H^3_{\et}(\X_K,\Z(2))$ as well as from $H^4_{\et}(\X_K,\Z(2))$. The examples
$H^1_{\et}(\X_K,\Z(1))\cong K^\times$, or $H^2_\et(\X_K,\Z(1))\cong \Pic(\X_K)$,
an extension of a finitely generated group by a finitely generated $\Z_p$-module, 
also suggest that the cohomology groups are topological groups. 
Thus our goal
is to construct topological cohomology groups which agree with \'etale 
cohomology groups with finite coefficients, but satisfy a Pontrjagin duality. 
More generally, we conjecture the existence of a cohomology theory on the category 
of separated schemes 
of finite type over $\mathrm{Spec}(\mathcal{O}_K)$, whose main expected properties 
are outlined in the last section of this paper. 

Its existence was suggested by the "Weil-Arakelov cohomology" of arithmetic schemes, 
which is conditionally defined in \cite{Flach-Morin-17} for proper regular schemes over 
$\mathrm{Spec}(\Z)$.
The aim of this paper is to give a possible construction of such groups.

Let $\mathrm{LCA}$ be the quasi-abelian category of locally compact abelian groups, 
and $\mathrm{FLCA}\subseteq \mathrm{LCA}$ the full subcategory  consisting of locally 
compact abelian groups of finite ranks in the sense of \cite{Hoffmann-Spitzweck-07}.
We consider the bounded derived category $\mathbf{D}^b(\mathrm{LCA})$ 
and $\mathbf{D}^b(\mathrm{FLCA})$, repectively, \cite{Hoffmann-Spitzweck-07}. 
The category $\mathbf{D}^b(\mathrm{FLCA})$ is a closed symmetric monoidal category
with internal homomorphisms  
$R\underline{\mathrm{Hom}}(-,-)$ and tensor product $\underline{\otimes}^{\bL}$.   
Assuming certain expected properties of Bloch's cycle complex $\Z(n)$, 
we construct for any $n\in \Z$ complexes in $\mathbf{D}^b(\mathrm{LCA})$
fitting in an exact triangle
\begin{equation}\label{IntroTri}
R\Gamma_{ar}(\X_s,Ri^!\Z(n))\rightarrow R\Gamma_{ar}(\X,\Z(n))\rightarrow R\Gamma_{ar}(\X_K,\Z(n)),
\end{equation}
and we define  
$$R\Gamma_{ar}(-,\br/\Z(n)):=R\Gamma_{ar}(-,\Z(n))\underline{\otimes}^{\bL} \br/\Z.$$ 
We expect that this theory satisfies duality and many other properties,
see Section 6.

To obtain unconditional results, we give an alternative construction, 
which conjecturally agrees with the above construction
of the triangle (\ref{IntroTri}) for $n=0,d$, and show that this triangle belongs to 
$\mathbf{D}^b(\mathrm{FLCA})$.
Then we prove Theorem \ref{dualityconj} below under the following hypothesis. 
 
\begin{Hypothesis}\label{conjfgintro}
The reduced scheme $(\X_s)^\red$ is a simple normal crossing  scheme, 
and the complex $R\Gamma_W(\mathcal{X}_s,\Z^c(0))$ is a perfect complex 
of abelian groups, where $\Z^c(0)$ denotes the Bloch cycle complex in 
its homological notation \cite{Geisser10} and $R\Gamma_W(\mathcal{X}_s,-)$ 
denotes Weil-\'etale cohomology.
\end{Hypothesis}

The homology groups of the complex $R\Gamma_W(\mathcal{X}_{s},\Z^c(0))[1]$ 
are and called arithmetic homology with compact support and
denoted by $H^c_i(\mathcal{X}_{s,\mathrm{ar}},\Z)$  in \cite{Geisser10b}. 
For $d\leq 2$, 
the complex $R\Gamma_W(\mathcal{X}_{s},\Z^c(0))$ is perfect, 
and perfectness in general follows from a special case of Parshin's conjecture 
\cite[Proposition 4.2]{Geisser10b}.

\begin{thm}\label{dualityconj} Suppose that either $d\leq 2$ or that $\X_s$ satisfies Hypothesis \ref{conjfgintro}.  Then there is a trace map $H^{2d}_{ar}(\X_K,\mathbb{R}/\Z(d))\rightarrow \mathbb{R}/\Z$ and an equivalence
$$R\Gamma_{ar}(\X_K,\mathbb{R}/\Z(n))\stackrel{\sim}{\longrightarrow} R\underline{\mathrm{Hom}}(R\Gamma_{ar}(\X_K,\Z(d-n)), \mathbb{R}/\Z[-2d])$$
in  $\mathbf{D}^b(\mathrm{FLCA})$, for $n=0,d$.
\end{thm}
Combining Theorem \ref{dualityconj} with a consequence of Sato's work \cite{Sato20}, we obtain 
\begin{cor}\label{cor-lca}  

Suppose that $\X/\mathcal{O}_K$ has good or strictly semi-stable reduction and suppose that $R\Gamma_W(\mathcal{X}_s,\Z^c(0))$ is a perfect complex of abelian groups.
Then there is a perfect pairing of locally compact abelian groups
$$H^{i}_{ar}(\X_K,\mathbb{R}/\Z(n))\times H^{2d-i}_{ar}(\X_K,\Z(d-n))\rightarrow H^{2d}_{ar}(\X_K,\mathbb{R}/\Z(d))\rightarrow \mathbb{R}/\Z$$
for $n=0,d$ and any $i\in\Z$.
\end{cor}

In a forthcoming paper we will give applications of this result to class
field theory of schemes over local fields.

\section{Locally compact abelian groups}\label{sectionLCA}

In this section we define and study the derived $\infty$-categories $\mathbf{D}^b(\mathrm{LCA})$ and $\mathbf{D}^b(\mathrm{FLCA})$. We also introduce a certain profinite completion functor. 

\subsection{Derived $\infty$-categories}

Let $A$ be an additive category. Let $\mathrm{C}^b(A)$ be the  differential graded 
category of bounded complexes of objects in $A$ and let 
$\mathcal{N}\subset \mathrm{C}^b(A)$ be a full subcategory which is closed under 
the formation of shifts and under the formation of mapping cones. 
If $\mathrm{N}_{\mathrm{dg}}(-)$ denotes the differential graded nerve 
\cite[Construction 1.3.1.6]{Lurie10}, then
$\mathrm{N}_{\mathrm{dg}}(\mathrm{C}^b(A))$ is a stable $\infty$-category and  
$\mathrm{N}_{\mathrm{dg}}(\mathcal{N})$ is a stable $\infty$-subcategory of 
$\mathrm{N}_{\mathrm{dg}}(\mathrm{C}^b(A))$ \cite[Proposition 1.3.2.10]{Lurie10}, 

The Verdier quotient is defined \cite[Theorem I.3.3]{Nikolaus-Scholze18} as the 
Dyer-Kan localization
$$\mathrm{N}_{\mathrm{dg}}(\mathrm{C}^b(A))/\mathrm{N}_{\mathrm{dg}}
(\mathcal{N}):=\mathrm{N}_{\mathrm{dg}}(\mathrm{C}^b(A))[W^{-1}],$$ 
where $W$ is the set of arrows in $\mathrm{N}_{\mathrm{dg}}(\mathrm{C}^b(A))$
whose cone lies in $\mathrm{N}_{\mathrm{dg}}(\mathcal{N})$.  
The $\infty$-category $\mathrm{N}_{\mathrm{dg}}(\mathrm{C}^b(A))[W^{-1}]$ is stable.
Moreover, the functor 
$\mathrm{N}_{\mathrm{dg}}(\mathrm{C}^b(A))\rightarrow \mathrm{N}_{\mathrm{dg}}
(\mathrm{C}^b(A))/\mathrm{N}_{\mathrm{dg}}(\mathcal{N})$ 
is exact and induces an equivalence from the category of exact functors 
$\mathrm{N}_{\mathrm{dg}}
(\mathrm{C}^b(A))/\mathrm{N}_{\mathrm{dg}}(\mathcal{N})\rightarrow\mathcal{E}$ 
to the category of exact functors 
$\mathrm{N}_{\mathrm{dg}}(\mathrm{C}^b(A))\rightarrow\mathcal{E}$ 
which send all objects of $\mathrm{N}_{\mathrm{dg}}(\mathcal{N})$ 
to zero objects in $\mathcal{E}$, 
for any (small) stable $\infty$-category $\mathcal{E}$. 
Finally, we have an equivalence of categories
$$h(\mathrm{N}_{\mathrm{dg}}(\mathrm{C}^b(A))/\mathrm{N}_{\mathrm{dg}}(\mathcal{N}))
\simeq h(\mathrm{N}_{\mathrm{dg}}(\mathrm{C}^b(A)))/h(\mathrm{N}_{\mathrm{dg}}
(\mathcal{N})))$$
where $h(-)$ denotes the homotopy category, 
and the right hand side is the classical Verdier quotient. 
Note that the homotopy category of a stable $\infty$-category is 
triangulated \cite[Theorem 1.1.2.14]{Lurie10}.

If $A$ is a quasi-abelian category in the sense of \cite{Schneiders99}, we define its bounded derived $\infty$-category
$$\mathbf{D}^b(A):=\mathrm{N}_{\mathrm{dg}}(\mathrm{C}^b(A))/\mathrm{N}_{\mathrm{dg}}(\mathcal{N})\simeq \mathrm{N}_{\mathrm{dg}}(\mathrm{C}^b(A))[S^{-1}]$$
where $\mathcal{N}\subset \mathrm{C}^b(A)$ is the full subcategory of strictly acyclic complexes, and $S$ is the set of strict quasi-isomorphisms. The homotopy category 
$$\mathrm{D}^b(A):= h(\mathbf{D}^b(A))$$
is equivalent to the bounded derived category of the quasi-abelian category $A$ in the sense of \cite{Schneiders99}.

\subsection{The category $\mathbf{D}^b(\mathrm{LCA})$.}\label{sectionDb(LCA)}
 We denote by $\mathrm{LCA}$ be the quasi-abelian category of locally compact abelian groups. A morphism of locally compact abelian groups $f:A\rightarrow B$ has a kernel $\mathrm{Ker}(f)=f^{-1}(0)$ and a cokernel $\mathrm{Coker}(f)=B/\overline{f(A)}$ where $\overline{f(A)}$ is the closure of $f(A)$ in $B$. The morphism $f$ is said to be strict if the map $\mathrm{Coker}(\mathrm{Ker}(f))\rightarrow \mathrm{Ker}(\mathrm{Coker}(f))$ is an isomorphism in $\mathrm{LCA}$. Then $f$ is strict if and only if the induced monomorphism $\overline{f}:A/\mathrm{Ker}(f)\rightarrow B$ is a closed embedding. Let  $\mathrm{FLCA}\subset \mathrm{LCA}$ be the quasi-abelian  category \cite[Corollary 2.11]{Hoffmann-Spitzweck-07} of locally compact abelian groups of finite ranks in the sense of \cite[Definition 2.6]{Hoffmann-Spitzweck-07}. Recall that $A\in \mathrm{LCA}$ has finite ranks if the $\br$-vector spaces of continuous morphisms $\underline{\mathrm{Hom}}(\br,A)$ and $\underline{\mathrm{Hom}}(A,\br)$ are finite dimensional and $p:A\rightarrow A$ is strict with finite kernel and cokernel for any prime number $p$. 
 
Let $\mathbf{D}^b(\mathrm{LCA})$ and $\mathbf{D}^b(\mathrm{FLCA})$ be the bounded derived $\infty$-category of $\mathrm{LCA}$ and $\mathrm{FLCA}$, respectively. 
Then $\mathbf{D}^b(\mathrm{LCA})$ and $\mathbf{D}^b(\mathrm{FLCA})$ are stable 
$\infty$-categories in the sense of \cite{Lurie10} whose homotopy categories are the 
bounded derived categories $\mathrm{D}^b(\mathrm{LCA})$ and 
$\mathrm{D}^b(\mathrm{FLCA})$  as defined in \cite{Hoffmann-Spitzweck-07}, respectively. It is more convenient to work with the derived $\infty$-category $\mathbf{D}^b(\mathrm{LCA})$ rather than with its homotopy category. For example, let $\mathrm{Fun}(\Delta^1,\mathbf{D}^b(\mathrm{LCA}))$ be the $\infty$-category of arrows in $\mathbf{D}^b(\mathrm{LCA})$. Taking the mapping fiber (or cofiber) of a morphism defines a functor (see \cite[Remark 1.1.1.7]{Lurie10})
$$\fonc{\mathrm{Fib}}{\mathrm{Fun}(\Delta^1,\mathbf{D}^b(\mathrm{LCA}))}{\mathbf{D}^b(\mathrm{LCA})}{C\rightarrow C'}{C\times_{C'}0}.
$$ 
Let $\mathrm{TA}$ be the quasi-abelian category of topological abelian groups, and define $\mathbf{D}^b(\mathrm{TA})$ and $\mathrm{D}^b(\mathrm{TA})$ as above.
The inclusions $\mathrm{FLCA}\subset \mathrm{LCA}\subset\mathrm{TA}$ send strict quasi-isomorphisms to strict quasi-isomorphisms, hence induce functors
$$\mathbf{D}^b(\mathrm{FLCA})\rightarrow \mathbf{D}^b(\mathrm{LCA})\rightarrow \mathbf{D}^b(\mathrm{TA}).$$ 
The functor $\mathrm{disc}:\mathrm{TA}\rightarrow \mathrm{Ab}$, sending a topological abelian group to its underlying discrete abelian group, sends strict quasi-isomorphisms to usual quasi-isomorphisms. This yields a functor
$$\mathrm{disc}:\mathbf{D}^b(\mathrm{TA})\rightarrow \mathbf{D}^b(\mathrm{Ab}).$$
Recall that the Pontryagin dual $X^D:=\underline{\mathrm{Hom}}(X,\br/\Z)$ of the locally compact abelian group $X$ is the group of continuous homomorphisms $X\rightarrow\br/\Z$ endowed with the compact-open topology, and that Pontryagin duality gives an isomorphism of locally compact groups
$$X\stackrel{\sim}{\rightarrow} X^{DD}.$$
The functor $(-)^D$ sends strict quasi-isomorphisms to strict quasi-isomorphisms and locally compact compact abelian groups of finite ranks to locally compact groups of  finite ranks. We obtain equivalences of $\infty$-categories
$$
\appl{\mathbf{D}^b(\mathrm{LCA})^{\mathrm{op}}}{\mathbf{D}^b(\mathrm{LCA})}{X}{X^D}
$$
and
$$
\appl{\mathbf{D}^b(\mathrm{FLCA})^{\mathrm{op}}}{\mathbf{D}^b(\mathrm{FLCA})}{X}{X^D}.
$$
In \cite{Hoffmann-Spitzweck-07}, the authors define functors
$$R\mathrm{Hom}_{\mathrm{LCA}}(-,-):\mathrm{D}^b(\mathrm{LCA})^{\mathrm{op}}\times \mathrm{D}^b(\mathrm{LCA})\rightarrow \mathrm{D}^b(\mathrm{TA})$$
and
$$R\mathrm{Hom}_{\mathrm{FLCA}}(-,-):\mathrm{D}^b(\mathrm{FLCA})^{\mathrm{op}}\times \mathrm{D}^b(\mathrm{FLCA})\rightarrow \mathrm{D}^b(\mathrm{FLCA}).$$
The construction of the functor $R\mathrm{Hom}_{\mathrm{FLCA}}(-,-)$ actually gives a functor of stable $\infty$-categories
$$R\underline{\mathrm{Hom}}(-,-):\mathbf{D}^b(\mathrm{FLCA})^{\mathrm{op}}\times \mathbf{D}^b(\mathrm{FLCA})\rightarrow \mathbf{D}^b(\mathrm{FLCA}).$$
Indeed, let $\mathrm{I}$ (resp. $\mathrm{P}$) be the additive category of divisible (resp. codivisible) locally compact abelian groups $I$ (resp. $P$) of finite ranks such that $I_{\Z}=0$ (such that $P_{\mathbb{S}^1}=0$), see \cite[Definition 3.2]{Hoffmann-Spitzweck-07}.
Define 
$$\mathbf{D}^b(\mathrm{I}):= \mathrm{N}_{\mathrm{dg}}(\mathrm{C}^b(\mathrm{I}))/\mathrm{N}_{\mathrm{dg}}(\mathcal{N}_{\mathrm{I}})$$
where $\mathcal{N}_{\mathrm{I}}\subset \mathrm{C}^b(\mathrm{I})$ is the $dg$-subcategory of strictly acyclic bounded complexes. We define similarly
$$\mathbf{D}^b(\mathrm{P}):= \mathrm{N}_{\mathrm{dg}}(\mathrm{C}^b(\mathrm{P}))/\mathrm{N}_{\mathrm{dg}}(\mathcal{N}_{\mathrm{P}}).$$
The exact functor
 $$\mathrm{N}_{\mathrm{dg}}(\mathrm{C}^b(\mathrm{I}))\rightarrow \mathrm{N}_{\mathrm{dg}}(\mathrm{C}^b(\mathrm{FLCA}))\rightarrow \mathbf{D}^b(\mathrm{FLCA})$$
 induces an exact functor
 \begin{equation}\label{feq}
\mathbf{D}^b(\mathrm{I})\rightarrow \mathbf{D}^b(\mathrm{FLCA})
\end{equation}
of stable $\infty$-categories which induces an equivalences between the corresponding homotopy categories by \cite[Cor. 3.10]{Hoffmann-Spitzweck-07}. It follows that (\ref{feq}) is an equivalence of stable $\infty$-categories. Similarly
$\mathbf{D}^b(\mathrm{P})\rightarrow \mathbf{D}^b(\mathrm{FLCA})$
is an equivalence. We may therefore define
$$R\underline{\mathrm{Hom}}(-,-):\mathbf{D}^b(\mathrm{FLCA})^{\mathrm{op}}\times\mathbf{D}^b(\mathrm{FLCA})\stackrel{\sim}{\leftarrow}\mathbf{D}^b(\mathrm{P})^{\mathrm{op}}\times\mathbf{D}^b(\mathrm{I})\rightarrow \mathbf{D}^b(\mathrm{FLCA})$$
since the functor
$$
\appl{\mathrm{C}^b(\mathrm{P})^{\mathrm{op}}\times\mathrm{C}^b(\mathrm{I})}{ \mathrm{C}^b(\mathrm{FLCA})}{(P,I)}{\underline{\mathrm{Hom}}^{\bullet}(P,I):=\mathrm{Tot}(\underline{\mathrm{Hom}}(P,I))}
$$
sends a pair of strict quasi-isomorphisms to a strict quasi-isomorphism \cite[Cor. 3.7]{Hoffmann-Spitzweck-07}. Here $\underline{\mathrm{Hom}}(P,I)$ is the double complex of continuous maps endowed with the compact-open topology, and $\mathrm{Tot}$ denotes the total complex. Note that the Pontryagin dual $X^D$ is given by the functor
$$\fonc{R\underline{\mathrm{Hom}}(-,\mathbb{R}/\Z)}{\mathbf{D}^b(\mathrm{FLCA})^{\mathrm{op}}}{\mathbf{D}^b(\mathrm{FLCA})}{X}{X^D}.
$$
Following \cite{Hoffmann-Spitzweck-07}, we define the derived topological tensor product
\begin{equation}\label{tensorFLCA}
\appl{\mathbf{D}^b(\mathrm{FLCA})\times \mathbf{D}^b(\mathrm{FLCA})}{\mathbf{D}^b(\mathrm{FLCA})}{(X,Y)}{X\underline{\otimes}^{\bL}Y:= R\underline{\mathrm{Hom}}(X,Y^D)^D}.
\end{equation}

\begin{lem}\label{ff}
The functor $\mathbf{D}^b(\mathrm{FLCA})\rightarrow \mathbf{D}^b(\mathrm{LCA})$ is an exact and fully faithful functor of stable $\infty$-categories. 
\end{lem}
\begin{proof}
The functor $$\mathrm{N}_{\mathrm{dg}}(\mathrm{C}^b(\mathrm{FLCA}))\rightarrow\mathrm{N}_{\mathrm{dg}}(\mathrm{C}^b(\mathrm{LCA}))\rightarrow \mathbf{D}^b(\mathrm{LCA})$$
induces an exact functor $\mathbf{D}^b(\mathrm{FLCA})\rightarrow \mathbf{D}^b(\mathrm{LCA})$ by \cite[Theorem I.3.3(i)]{Nikolaus-Scholze18}. It remains to check that this functor is fully faithful.
The functors $R\mathrm{Hom}_{\mathrm{LCA}}(-,-)$ and $R\mathrm{Hom}_{\mathrm{FLCA}}(-,-)$ induce the same functor
$$\mathrm{D}^b(\mathrm{FLCA})^{\mathrm{op}}\times \mathrm{D}^b(\mathrm{FLCA})\rightarrow \mathrm{D}^b(\mathrm{TA})$$
by \cite[Remark 4.9]{Hoffmann-Spitzweck-07}. Moreover, for any $X,Y\in \mathrm{D}^b(\mathrm{FLCA})$ we have \cite[Proposition 4.12(i)]{Hoffmann-Spitzweck-07}
$$\mathrm{disc}(H^0(R\mathrm{Hom}_{\mathrm{LCA}}(X,Y)))\simeq \mathrm{Hom}_{\mathrm{D}^b(\mathrm{LCA})}(X,Y)$$
and similarly\footnote{One may adapt the proof of \cite[Proposition 4.12(i)]{Hoffmann-Spitzweck-07} to this case, using \cite[Corollary 3.10(iii)]{Hoffmann-Spitzweck-07}.}
$$\mathrm{disc}(H^0(R\mathrm{Hom}_{\mathrm{FLCA}}(X,Y)))\simeq \mathrm{Hom}_{\mathrm{D}^b(\mathrm{FLCA})}(X,Y).$$
Therefore, the map
$$\mathrm{Hom}_{\mathrm{D}^b(\mathrm{FLCA})}(X,Y)\rightarrow \mathrm{Hom}_{\mathrm{D}^b(\mathrm{LCA})}(X,Y)$$
is an isomorphism of abelian groups, i.e. $\mathrm{D}^b(\mathrm{FLCA})\rightarrow \mathrm{D}^b(\mathrm{LCA})$ is fully faithful. Hence
\begin{equation}\label{ff0}
\mathbf{D}^b(\mathrm{FLCA})\rightarrow \mathbf{D}^b(\mathrm{LCA})
\end{equation}
is an exact functor of stable $\infty$-categories which induces a fully faithful functor between the corresponding homotopy categories. It follows that (\ref{ff0}) is fully faithful.

\end{proof}
Therefore we may identify $\mathbf{D}^b(\mathrm{FLCA})$ with its 
essential image in  $\mathbf{D}^b(\mathrm{LCA})$.  
The stable $\infty$-category $\mathbf{D}^b(\mathrm{LCA})$ is endowed with 
a $t$-structure by \cite[Section 1.2.2]{Schneiders99}, 
since a $t$-structure on a stable $\infty$-category is defined as a 
$t$-structure on its homotopy category \cite[Definition 1.2.1.4]{Lurie10}.  
We denote its heart by $\mathcal{LH}(\mathrm{LCA})$. 
It is an abelian category containing $\mathrm{LCA}$ as a full subcategory
This also applies to $\mathbf{D}^b(\mathrm{FLCA})$, and we denote 
$\mathcal{LH}(\mathrm{FLCA})$ the heart of the corresponding $t$-structure.

\begin{rem}
By \cite[Cor.\ 1.2.21]{Schneiders99}, an object in 
$\mathcal{LH}(\mathrm{LCA})$ can be represented by a monomorphism 
$E_1\to E_0$ in $\mathrm{LCA}$, where $E_0$ is in degree zero. A common
example appearing below is a monomorphism $\Z^a\to \Z_p$. 
Its cokernel in $\mathrm{LCA}$ is trivial,
but for the underlying discrete abelian groups the cokernel is
$(\mathbb{Q}/\Z')\oplus (\mathbb{Q}/\Z)^{a-1}\oplus D$ 
with $D$ uniquely divisible.
\end{rem}

\begin{rem}\label{texact}
It follows from \cite[Corollary 2.11]{Hoffmann-Spitzweck-07} and  \cite[Proposition 1.2.19]{Schneiders99} that the fully faithful functor $\mathbf{D}^b(\mathrm{FLCA})\rightarrow\mathbf{D}^b(\mathrm{LCA})$ is $t$-exact.  Therefore, the induced functor $\mathcal{LH}(\mathrm{FLCA})\hookrightarrow\mathcal{LH}(\mathrm{LCA})$ is exact and fully faithful.
\end{rem}

\begin{notation}
For any $X\in \mathbf{D}^b(\mathrm{LCA})$ and any $i\in\Z$, we consider $$H^{i}(X):=\tau^{\geq 0}\tau^{\leq 0} (X[i]) \in \mathcal{LH}(\mathrm{LCA}).$$
\end{notation}
In view of Remark \ref{texact}, we identify $\mathcal{LH}(\mathrm{FLCA})$ with a full subcategory of $\mathcal{LH}(\mathrm{LCA})$. 

\begin{lem}\label{lemcoh}
Let $X\in \mathbf{D}^b(\mathrm{LCA})$. Then $X\in \mathbf{D}^b(\mathrm{FLCA})$ if and only if  $H^{i}(X)\in \mathcal{LH}(\mathrm{FLCA})$ for any $i\in\Z$.
\end{lem}
\begin{proof}

If $X\rightarrow Y\rightarrow Z$ is a fiber sequence in $\mathbf{D}^b(\mathrm{LCA})$ such that $X,Z\in \mathbf{D}^b(\mathrm{FLCA})$, then $Y\in \mathbf{D}^b(\mathrm{FLCA})$. Indeed, a stable subcategory is closed under extensions. 
Let $X\in \mathbf{D}^b(\mathrm{LCA})$ such that $H^{i}(X)\in \mathcal{LH}( \mathrm{FLCA})$ for any $i$. Note that $H^{i}(X)=0$ for all but finitely many $i\in\Z$. Therefore, $X$ has a finite exhaustive filtration with $i$-graded piece $H^{i}(X)[-i]\in \mathbf{D}^b(\mathrm{FLCA})$, so that $X$ belongs to $\mathbf{D}^b(\mathrm{FLCA})$ by induction.

The converse follows from the fact that the inclusion functor $\mathbf{D}^b(\mathrm{FLCA})\rightarrow \mathbf{D}^b(\mathrm{LCA})$ is $t$-exact by Remark \ref{texact}.
\end{proof}

Recall that Pontryagin duality gives an equivalence
$$
\appl{\mathbf{D}^b(\mathrm{LCA})^{\mathrm{op}}}{\mathbf{D}^b(\mathrm{LCA})}{X}{X^D}.
$$
\begin{lem}\label{remcohdual}
Let $X\in \mathbf{D}^b(\mathrm{LCA})$ such that
$H^{i}(X)\in \mathcal{LH}(\mathrm{LCA})$ belongs to $\mathrm{LCA}$ for any $i\in\Z$.
Then for any $i\in\Z$ we have a canonical isomorphism in $\mathrm{LCA}$
$$H^{i}(X^D)\simeq (H^{-i}(X))^D.$$
\end{lem}
\begin{proof}
Let $$X= [\cdots\rightarrow X^{i-1}\stackrel{d_X^{i-1}}{\rightarrow} X^{i}\stackrel{d_X^{i}}{\rightarrow} X^{i+1}\rightarrow \cdots ]$$
be an object of $\mathbf{D}^b(\mathrm{LCA})$ such that $H^{i}(X)\in \mathcal{LH}(\mathrm{LCA})$ belongs to $\mathrm{LCA}$ for any $i\in\Z$. We first observe that the differentials $d_X^{i}$ all are strict morphisms. By \cite[Proposition 1.2.19]{Schneiders99}, the object $H^{i}(X)$ of $\mathcal{LH}(\mathrm{LCA})$ is given by the complex $[0\rightarrow \mathrm{Coim}(d_X^{i-1})\stackrel{\delta}{\rightarrow} \mathrm{Ker}(d_X^{i})\rightarrow 0]$, where $\mathrm{Ker}(d_X^{i})$ sits in degree $0$ and $\delta$ is a monomorphism. Since $H^{i}(X)\in\mathrm{LCA}$, the map $\delta$ is strict by \cite[Proposition 1.2.29]{Schneiders99}, i.e. $\delta$ is a closed embedding. Since $\mathrm{Ker}(d_X^{i})\hookrightarrow X^{i}$ is a closed embedding as well, so is the map $\mathrm{Coim}(d_X^{i-1})=X^{i-1}/\mathrm{Ker}(d_X^{i-1})\rightarrow X^{i}$. Hence $d_X^{i-1}$ is strict.

We set $Y:=X^D$ so that $Y^{-i}=(X^{i})^D$ and $d_Y^{-i}:Y^{-i}\rightarrow Y^{-i+1}$ is the map $d_Y^{-i}=(d_X^{i-1})^D$. The differentials $d_X^{*}$ are all strict morphisms, hence so are their duals $d_Y^{*}$. We have the following isomorphisms of locally compact abelian groups:
\begin{eqnarray}
H^{i}(X)^D&\simeq&\left(\mathrm{Coker}(\mathrm{Coim}(d_X^{i-1})\stackrel{\delta}{\rightarrow} \mathrm{Ker}(d_X^{i}))\right)^D\label{dual1}\\
&\simeq&\mathrm{Ker}\left(\mathrm{Coker}((d_X^{i})^D)\rightarrow \mathrm{Im}((d_X^{i-1})^D)\right)\label{dual2}\\
&\simeq&\mathrm{Ker}\left(\mathrm{Coker}(d_Y^{-i-1})\rightarrow \mathrm{Im}(d_Y^{-i})\right)\label{dual3}\\
&\simeq&\mathrm{Ker}\left(\mathrm{Coker}(d_Y^{-i-1})\rightarrow Y^{-i+1}\right)\label{dual4}\\
&\simeq&\mathrm{Ker}\left(Y^{-i}/d_Y^{-i-1}(Y^{-i-1})\rightarrow Y^{-i+1}\right)\label{dual5}\\
&\simeq&\mathrm{Ker}(d_Y^{-i})/d_Y^{-i-1}(Y^{-i-1})\label{dual6}\\
&\simeq&\mathrm{Ker}(d_Y^{-i})/\mathrm{Coim}(d_Y^{-i-1})\label{dual7}\\
&\simeq& H^{-i}(Y) \label{dual8}
\end{eqnarray}
where the kernels, cokernels, images, and coimages are all computed in $\mathrm{LCA}$. The isomorphism (\ref{dual1}) is valid by \cite[Proposition 1.2.29]{Schneiders99} since the map $\delta$ is strict, and (\ref{dual2}) holds since Pontryagin duality $(-)^D:\mathrm{LCA}^{\mathrm{op}}\rightarrow \mathrm{LCA}$ is an equivalence of categories with kernels and cokernels. The identification (\ref{dual3}) is given by definition of the maps $d_Y^*$ and (\ref{dual4}) holds since $\mathrm{Im}(d_Y^{-i})\rightarrow Y^{-i+1}$ is a monomorphism. We have (\ref{dual5}) in view of $\mathrm{Coker}(d_Y^{-i-1})=Y^{-i}/d_Y^{-i-1}(Y^{-i-1})$, which is valid since $d_Y^{-i-1}(Y^{-i-1})$ is closed in  $Y^{-i}$, as $d_Y^{-i-1}$ is strict. The isomorphism of locally compact abelian groups (\ref{dual6}) is clear; (\ref{dual7}) holds since $\mathrm{Coim}(d_Y^{-i-1})\rightarrow d_Y^{-i-1}(Y^{-i-1})= \mathrm{Im}(d_Y^{-i-1})$ is an isomorphism in $\mathrm{LCA}$ since $d_Y^{-i-1}$ is strict. Finally, (\ref{dual8}) holds by \cite[Propositions 1.2.19 and 1.2.29]{Schneiders99} since $\mathrm{Coim}(d_Y^{-i-1})\rightarrow \mathrm{Ker}(d_Y^{-i})$ is strict; indeed $Y^{-i-1}/\mathrm{Ker}(d_Y^{-i-1})\rightarrow Y^{-i}$ is a closed embedding hence so is $Y^{-i-1}/\mathrm{Ker}(d_Y^{-i-1})\rightarrow \mathrm{Ker}(d_Y^{-i})$.
\end{proof}

$\mathcal{LH}(\mathrm{LCA})^{\mathrm{op}}
\stackrel{\sim}{\rightarrow}\mathcal{RH}(\mathrm{LCA})$, where 
$\mathcal{RH}(\mathrm{LCA})$ is the heart of the right $t$-structure on 
$\mathbf{D}^b(\mathrm{LCA})$ which we do not consider in this paper.

The inclusion $\mathrm{Ab}\subset \mathrm{LCA}$ as discrete objects
induces an exact functor
$$i:\mathbf{D}^b(\mathrm{Ab})\rightarrow \mathbf{D}^b(\mathrm{LCA}).$$

\begin{prop}
The exact functor $i:\mathbf{D}^b(\mathrm{Ab})\rightarrow \mathbf{D}^b(\mathrm{LCA})$ is fully faithful and left adjoint to 
$$\mathrm{disc}:\mathbf{D}^b(\mathrm{LCA})\rightarrow \mathbf{D}^b(\mathrm{Ab}).$$
\end{prop}

\begin{proof}
The functor
$$\mathrm{C}^b(\mathrm{Ab})\stackrel{i}{\rightarrow} \mathrm{C}^b(\mathrm{LCA})\stackrel{\mathrm{disc}}{\rightarrow} \mathrm{C}^b(\mathrm{Ab})$$
is isomorphic to the identity functor of $\mathrm{C}^b(\mathrm{Ab})$. We obtain a natural transformation
\begin{equation}\label{unit}
\mathrm{Id}_{\mathbf{D}^b(\mathrm{Ab})}\stackrel{\sim}{\rightarrow} \mathrm{disc}\circ i .
\end{equation}
Similarly, there is a natural transformation
$$i\circ \mathrm{disc}\rightarrow \mathrm{Id}_{\mathbf{D}^b(\mathrm{LCA})}.$$

Let $X\in \mathbf{D}^b(\mathrm{Ab})$ and let $Y\in \mathbf{D}^b(\mathrm{LCA})$. Let $F\stackrel{\sim}{\rightarrow}X$ be a bounded flat resolution, and let $Y\stackrel{\sim}{\rightarrow} D$ be a strict quasi-isomorphism where $D$ is a bounded complex of divisible locally compact abelian groups. Then $F$ is a bounded complex of codivisible\footnote{$A\in\mathrm{LCA}$ is said to be codivisible if $A^D$ is divisible.}  discrete groups $F^{i}$ (in particular, such that $F^{i}_{\mathbb{S}^1}=0$). Therefore, we have
$$R\mathrm{Hom}_{\mathrm{LCA}}(i(X),Y)\simeq \underline{\mathrm{Hom}}^{\bullet}(F,D):=\mathrm{Tot}(\underline{\mathrm{Hom}}(F,D))$$
by \cite[Corollary 4.7]{Hoffmann-Spitzweck-07}, where $\underline{\mathrm{Hom}}(F,D)$ is the double complex of continuous maps endowed with the compact-open topology, and $\mathrm{Tot}$ denotes the total complex.
We obtain 
\begin{eqnarray*}
\mathrm{disc}\left(R\mathrm{Hom}_{\mathrm{LCA}}(i(X),Y)\right)&\simeq & \mathrm{disc}\left(\underline{\mathrm{Hom}}^{\bullet}(F,D)\right)\\
&\simeq& \mathrm{Hom}^{\bullet}(F,\mathrm{disc}(D))\\
&\simeq& R\mathrm{Hom}(X,\mathrm{disc}(Y)).
\end{eqnarray*}
In view of \cite[Proposition 4.12]{Hoffmann-Spitzweck-07} we have
\begin{eqnarray*}
H^0(\mathrm{disc}\left(R\mathrm{Hom}_{\mathrm{LCA}}(i(X),Y[-n])\right))&\simeq& \mathrm{disc}(H^0\left(R\mathrm{Hom}_{\mathrm{LCA}}(i(X),Y[-n])\right))\\
&\simeq&\mathrm{Hom}_{\mathrm{D}^b(\mathrm{LCA})}(i(X),Y[-n])\\
&\simeq&\pi_0(\mathrm{Map}_{\mathbf{D}^b(\mathrm{LCA})}(i(X),\Omega^nY))\\
&\simeq& \pi_n( \mathrm{Map}_{\mathbf{D}^b(\mathrm{LCA})}(i(X),Y))
\end{eqnarray*}
where $\Omega(-):=0\times_{(-)}0$ is the loop space functor. Similarly we have
$$H^0(R\mathrm{Hom}(X,\mathrm{disc}(Y[-n])))\simeq \pi_n( \mathrm{Map}_{\mathbf{D}^b(\mathrm{Ab})}(X,\mathrm{disc}(Y))).$$
Hence the map
$$\mathrm{Map}_{\mathbf{D}^b(\mathrm{LCA})}(i(X),Y)\rightarrow \mathrm{Map}_{\mathbf{D}^b(\mathrm{Ab})}(X,\mathrm{disc}(Y))$$
is an equivalence of $\infty$-groupoids. The result then follows from \cite[Proposition 5.2.2.8]{HTT} and from the fact that the unit transformation (\ref{unit}) is an equivalence.
\end{proof}

\begin{defn}
An object $X\in \mathbf{D}^b(\mathrm{LCA})$ lies in the essential image of the functor  $i:\mathbf{D}^b(\mathrm{Ab})\rightarrow \mathbf{D}^b(\mathrm{LCA})$ if and only if the co-unit map
$i\circ \mathrm{disc}(X)\rightarrow X$ is an equivalence. Such an object $X\in \mathbf{D}^b(\mathrm{LCA})$ is called \emph{discrete}.
\end{defn}

\begin{lem}\label{lemmapdiscretecont}
Let $X,Y\in \mathbf{D}^b(\mathrm{Ab})$. If $iX$ and $iY$ belong to $
\mathbf{D}^b(\mathrm{FLCA})$, then there 
is a canonical map
$$i(R\mathrm{Hom}(X,Y))\rightarrow R\underline{\mathrm{Hom}}(iX,iY).$$
Moreover, if $X,Y$ are perfect complexes of abelian groups, then this map is an equivalence.
\end{lem}
\begin{proof}
Let $P\stackrel{\sim}{\rightarrow} iX$ and $iY\stackrel{\sim}{\rightarrow} I$ be strict quasi-isomorphisms where $P\in\mathrm{C}^b(\mathrm{P})$ (resp. $I\in\mathrm{C}^b(\mathrm{I})$). We denote by $P^{\delta}:=\mathrm{disc}(P)$ and $I^{\delta}:=\mathrm{disc}(P)$ the underlying complexes of discrete abelian groups. Then the maps $P^{\delta}\stackrel{\sim}{\rightarrow} X$ and $Y\stackrel{\sim}{\rightarrow} I^{\delta}$ are quasi-isomorphisms in the usual sense. Hence we have
$\mathrm{Hom}^{\bullet}(X,I^{\delta})\simeq R\mathrm{Hom}(X,Y)$, where $\mathrm{Hom}^{\bullet}$ denotes the total complex of the double complex of morphisms of abelian groups. We denote by $\underline{\mathrm{Hom}}^{\bullet}(P,I)$ the total complex of the double complex of continuous morphisms endowed with the compact-open topology. Then we have morphisms
$$R\mathrm{Hom}(X,Y)\simeq \mathrm{Hom}^{\bullet}(X,I^{\delta})\rightarrow \underline{\mathrm{Hom}}^{\bullet}(iX,I)\rightarrow \underline{\mathrm{Hom}}^{\bullet}(P,I)\simeq R\underline{\mathrm{Hom}}(iX,iY).$$

Suppose now that $X$ and $Y$ are perfect complexes of abelian groups. We may suppose that $X^{n}$ is a finitely generated free abelian group for all $n\in\Z$, zero for almost all $n$, and similarly for $Y$. 
We have a strict quasi-isomorphism $$iY\stackrel{\sim}{\rightarrow} I:=\mathrm{Tot}[Y\otimes\mathbb{R}\rightarrow Y\otimes\mathbb{R}/\Z]$$ 
where $[Y\otimes\mathbb{R}\rightarrow Y\otimes\mathbb{R}/\Z]$ is seen as a double complex of locally compact abelian groups and $\mathrm{Tot}$ is the total complex. Then $iX\in \mathrm{C}^b(\mathrm{P})$ and $I\in \mathrm{C}^b(\mathrm{I})$ and we have a strict quasi-isomorphism
$$i\mathrm{Hom}^{\bullet}(X,Y)\stackrel{\sim}{\rightarrow} \underline{\mathrm{Hom}}^{\bullet}(iX,I).$$
We obtain
$$iR\mathrm{Hom}(X,Y)\simeq i\mathrm{Hom}^{\bullet}(X,Y)\stackrel{\sim}{\rightarrow} \underline{\mathrm{Hom}}^{\bullet}(iX,I)\simeq R\underline{\mathrm{Hom}}(iX,iY).$$

\end{proof}

\subsection{Profinite completion}

\begin{defn} We define a functor
$$\fonc{(-)\widehat{\underline{\otimes}}\widehat{\Z}}{\mathbf{D}^b(\mathrm{Ab})}{\mathbf{D}^b(\mathrm{LCA})}{X}{\left(i(\mathrm{colim}\,R\mathrm{Hom}(X,\Z/m)\right))^D},$$
where we  compute $R\mathrm{Hom}(X,\Z/m)$ and the colimit $\mathrm{colim}\,R\mathrm{Hom}(X,\Z/m)$ over $m$ in the $\infty$-category $ \mathbf{D}^b(\mathrm{Ab})$. 
We define similarly
$$\fonc{(-)\widehat{\underline{\otimes}}\Z_p}{\mathbf{D}^b(\mathrm{Ab})}{\mathbf{D}^b(\mathrm{LCA})}{X}{\left(i(\mathrm{colim}\,R\mathrm{Hom}(X,\Z/p^{\bullet})\right))^D}.$$
\end{defn}
For any $X\in \mathbf{D}^b(\mathrm{LCA})$ we define\footnote{This is compatible with the definition given in Section \ref{sectionDb(LCA)}, which is only valid if $X\in \mathbf{D}^b(\mathrm{FLCA})$.}
$$R\underline{\mathrm{Hom}}(X,\Z/m):=\mathrm{Fib}(X^D\stackrel{m}{\longrightarrow} X^D)$$
and
$$X\underline{\otimes}^{\bL}\Z/m:=\mathrm{Cofib}(X\stackrel{m}{\longrightarrow} X).$$

\begin{prop}\label{discprocompletion}
Let $X\in \mathbf{D}^b(\mathrm{Ab})$. Suppose that $R\underline{\mathrm{Hom}}(i(X),\Z/m)\in \mathbf{D}^b(\mathrm{LCA})$ is discrete for any $m$. Then we have an equivalence
$$X\widehat{\underline{\otimes}}\widehat{\Z}\simeq \underleftarrow{\mathrm{lim}}\left(i(X)\underline{\otimes}^{\bL}\Z/m\right)$$
where the limit is computed in the $\infty$-category $\mathbf{D}^b(\mathrm{LCA})$ and an equivalence
$$\mathrm{disc}(X\widehat{\underline{\otimes}}\widehat{\Z})\simeq X\widehat{\otimes}\widehat{\Z}:=\underleftarrow{\mathrm{lim}}(X\otimes^{\bL}\Z/m)\in \mathbf{D}^b(\mathrm{Ab}).$$
\end{prop}
\begin{proof}
The co-unit map
$$i\circ \mathrm{disc} \,R\underline{\mathrm{Hom}}(i(X),\Z/m)\rightarrow R\underline{\mathrm{Hom}}(i(X),\Z/m)$$
is an equivalence by assumption. Moreover we have 
$$R\mathrm{Hom}(X,\Z/m) \simeq\mathrm{disc}\, R\underline{\mathrm{Hom}}(i(X),\Z/m)$$
hence
$$i\,R\mathrm{Hom}(X,\Z/m)\stackrel{\sim}{\rightarrow} R\underline{\mathrm{Hom}}(i(X),\Z/m).$$
We obtain
\begin{eqnarray*}
X\widehat{\underline{\otimes}}\widehat{\Z}&:=&\left(i(\mathrm{colim}\,R\mathrm{Hom}(X,\Z/m)\right))^D\\
&\simeq & \left(\mathrm{colim}\,(i\,R\mathrm{Hom}(X,\Z/m)\right))^D\\
&\simeq & \mathrm{lim} \left((i\,R\mathrm{Hom}(X,\Z/m))^D\right)\\
&\simeq & \mathrm{lim} \left(R\underline{\mathrm{Hom}}(i(X),\Z/m)^D\right)\\
&\simeq & \mathrm{lim} \left(i(X)\underline{\otimes}^{\bL} \Z/m\right)
\end{eqnarray*}
since the left adjoint functor $i$ commutes with arbitrary colimits, and since  $(-)^D$ transforms colimits into limits. Hence we have
\begin{eqnarray*}
\mathrm{disc}(X\widehat{\underline{\otimes}}\widehat{\Z}) &\simeq & \mathrm{disc}\left(\mathrm{lim} (i(X)\underline{\otimes}^{\bL} \Z/m) \right)\\
&\simeq & \mathrm{lim}\left(\mathrm{disc} (i(X)\underline{\otimes}^{\bL} \Z/m) \right)\\
&\simeq & \mathrm{lim}\left(\mathrm{Cofib}(\mathrm{disc} \circ i(X)\stackrel{m}{\longrightarrow} \mathrm{disc} \circ i(X)) \right)\\
&\simeq & \mathrm{lim}\left(X\otimes^{\bL} \Z/m \right)
\end{eqnarray*}
since the right adjoint functor $\mathrm{disc}$ commutes with arbitrary limits.
\end{proof}
\begin{rem}\label{remfordiscrete}
Suppose that $X\in \mathbf{D}^b(\mathrm{Ab})$ is such that the cohomology groups of $X\otimes^{\bL}\Z/m$ are all finite. Then $R\underline{\mathrm{Hom}}(i(X),\Z/m)$ is discrete.
\end{rem}

\begin{rem} \label{remdef2pro}
We have 
$$X\widehat{\underline{\otimes}}\widehat{\Z}\simeq R\underline{\mathrm{Hom}}\left(i\,\mathrm{colim}\,R\mathrm{Hom}(X\otimes^\bL\Z/m,\mathbb{Q}/\Z),\mathbb{R}/\Z\right).$$
\end{rem}
\begin{lem} \label{lemmaptohat}
We have a canonical map $iX\rightarrow X\widehat{\underline{\otimes}}\widehat{\Z}$ in $\mathbf{D}^b(\mathrm{LCA})$.
\end{lem}
\begin{proof}
The composite map
$$i(R\mathrm{Hom}(X,\Z/m))\stackrel{\sim}{\rightarrow}i\circ \mathrm{disc} (R\underline{\mathrm{Hom}}(iX,\Z/m))\rightarrow R\underline{\mathrm{Hom}}(iX,\Z/m))$$
$$\rightarrow R\underline{\mathrm{Hom}}(iX,\mathbb{R}/\Z)\simeq (iX)^D$$
induces
$$i(\mathrm{colim}\,R\mathrm{Hom}(X,\Z/m))\simeq \mathrm{colim}\,i(R\mathrm{Hom}(X,\Z/m))\rightarrow (iX)^D.$$
We obtain
$$iX\stackrel{\sim}{\rightarrow} (iX)^{DD}\rightarrow \left(i(\mathrm{colim}\,R\mathrm{Hom}(X,\Z/m))\right)^D=:X\widehat{\underline{\otimes}}\widehat{\Z}.$$
\end{proof}

\begin{rem}
Let $X$ be an object of $\mathbf{D}^b(\mathrm{Ab})$ whose image $iX\in\mathbf{D}^b(\mathrm{LCA})$ belongs to $\mathbf{D}^b(\mathrm{FLCA})$. Then one may consider $iX\underline{\otimes}^{\bL}\widehat{\Z}$ and $iX\underline{\otimes}^{\bL}\widehat{\Z}_p$ where $\underline{\otimes}^{\bL}$ is the tensor product (\ref{tensorFLCA}) in $\mathbf{D}^b(\mathrm{FLCA})$. There are canonical maps $iX\underline{\otimes}^{\bL}\widehat{\Z}\rightarrow X\widehat{\underline{\otimes}}\widehat{\Z}$ and $iX\underline{\otimes}^{\bL}\Z_p\rightarrow X\widehat{\underline{\otimes}}\Z_p$ 
but those maps are not equivalences in general. For example, we have
$$\mathbb{Q}_p/\Z_p \underline{\otimes}^{\bL} \Z_p\simeq \mathbb{Q}_p/\Z_p  \hspace{0.2cm}\textrm{ and } \hspace{0.2cm} \mathbb{Q} \underline{\otimes}^{\bL} \Z_p\simeq \mathbb{Q}_p$$
while
$$\mathbb{Q}_p/\Z_p \widehat{\underline{\otimes}} \Z_p\simeq \Z_p[1] \hspace{0.2cm}\textrm{ and } \hspace{0.2cm} \mathbb{Q} \widehat{\underline{\otimes}} \Z_p\simeq 0.$$ 
\end{rem}

\begin{notation}
In the next sections, given $X\in \mathbf{D}^b(\mathrm{Ab})$ we often simply denote by $X$ its image $iX$ in $\mathbf{D}^b(\mathrm{LCA})$. In particular, for $X,Y\in \mathbf{D}^b(\mathrm{Ab})$, we denote by $R\mathrm{Hom}(X,Y)\in \mathbf{D}^b(\mathrm{Ab})\subseteq \mathbf{D}^b(\mathrm{LCA})$ the usual $R\mathrm{Hom}$ seen as an object of $\mathbf{D}^b(\mathrm{LCA})$.
\end{notation}

\section{Duality for schemes over finite fields}

Let $Y$ be a proper scheme over a finite field. If $Y$ is smooth, then
the Weil-\'etale  cohomology $R\Gamma_{W}(Y,\Z)$ of \cite{Geisser06} 
is a perfect complex of abelian groups. In general, the Weil-h cohomology 
$R\Gamma_{Wh}(Y,\Z)$ of \cite{Geisser06} is a perfect complex of 
abelian groups provided resolution of singularities 
\cite[Definition 2.4]{Geisser06} holds 
(see Proposition \ref{propperfectWh} below). 
We show that if $Y$ is a simple normal crossing  scheme, then 
the Weil-\'etale cohomology $R\Gamma_{W}(Y,\Z)$ is a perfect complex 
of abelian groups, and that the canonical map 
$R\Gamma_{W}(Y,\Z)\rightarrow R\Gamma_{Wh}(Y,\Z)$ is 
an equivalence under resolution of singularities. 
In Section \ref{sectdualityXs}, we show that $R\Gamma_{W}(Y,\Z)$ is dual 
to $R\Gamma_{W}(Y,\Z^c(0))$ under the assumption that $R\Gamma_{W}(Y,\Z^c(0))$ 
is perfect, where $\Z^c(0)$ is the cycle complex.

\subsection{Finite generation of cohomology}
\label{sect-Wh-snc}

\begin{defn}\label{defnW*}
Let $k$ be a finite field with algebraic closure $\bar k$ and $W_k$ be its 
Weil-group. For a scheme $Y$ over $k$ we let $\bar Y=Y\times_k\bar k$.
For a scheme $Y$ of finite type $Y$ over $k$ we define the $Wh$-cohomology of 
the constant sheaf $\Z$ to be
$$R\Gamma_{Wh}(Y,\Z):= R\Gamma(W_k, R\Gamma_{eh}(\bar Y,\Z)).$$
\end{defn}

\begin{prop}\label{propperfectWh}
Let $Y$ be a proper scheme over a finite field $k$. 
Assume resolution of singularities for schemes over $k$ of dimension 
$\leq \mathrm{dim}(Y)$ \cite[Definition 2.4]{Geisser06}. 
Then $R\Gamma_{Wh}(Y,\Z)$ is a perfect complex of abelian groups.
\end{prop}

\begin{proof}
To prove perfectness of $R\Gamma_{Wh}(Y,\Z)$, one first reduces to the 
smooth and projective case by \cite[Prop. 3.2]{Geisser06}, in which case
one can conclude with loc. cit. Theorem 4.3 and \cite{Lichtenbaum05}.
\end{proof}

\begin{defn}\label{strictNCS}
Let $k$ be a field and let $Y$ be a pure dimensional proper scheme over $k$
with irreducible components $Y_i$, $i=1,\ldots,c$. 
Then $Y$ is said to be a \emph{simple normal crossing scheme} if 
for all  $I\subseteq \{1,\ldots, c\}$, 
$Y_I = \bigcap _{i\in I} Y_i$ is regular of codimension
$ |I|-1 $ in $Y$.
\end{defn} 

In fact for all the results in this paper we only need that
$(Y_I)^\red$ is regular.

\begin{lem}\label{finiteblowup}
Consider a cartesian square
$$\begin{CD}
Y'@>i'>> T'\\
@Vq VV@VV\pi V\\
Y@>i>> T
\end{CD}$$
of schemes of finite type over a field $k$ with $i$ a closed embedding
and $\pi$ finite such that $\pi|{T'-Y'}$ is an isomorphism to $T-Y$.
Then there is a distinguished triangle
$$R\Gamma_{et}(T,\Z)\to R\Gamma_{et}(T',\Z)\oplus 
R\Gamma_{et}(Y,\Z)\to R\Gamma_{et}(Y',\Z).$$
In particular, if $k$ is a finite field, we obtain a triangle
$$R\Gamma_{W}(T,\Z)\to R\Gamma_{W}(T',\Z)\oplus 
R\Gamma_{W}(Y,\Z)\to R\Gamma_{W}(Y',\Z).$$
\end{lem}

\begin{proof}
To get the first triangle, noting that $i_*$ and $\pi_*$ are exact,
it suffices to show that the exact sequence 
$$ 0\to \Z\to \pi_*\Z\oplus i_*\Z \to (\pi\circ i')_*\Z\to 0$$
of \'etale sheaves on $T$ is exact. But this follows by considering 
stalks at points $t\in T$. If $t\not\in Y$, then the sequence reduces
to the isomorphism $\Z\cong \pi_*\Z$, and if $t\in Y$, then 
$\Z\cong i_*\Z$ and $\pi_*\Z \cong (\pi\circ i')_*\Z$.

The second triangle can be obtain by applying $R\Gamma(W,-)$ to the 
first triangle after base extension to the algebraic closure.
\end{proof}

\begin{prop}\label{corW*}\label{prop-normalcrossing}
If $T^\red$ is a strict normal crossing scheme, then 
$R\Gamma_{W}(T,\Z)$ is a perfect complex of abelian groups. 
Under resolution of singularities, we have a quasi-isomorphic
$R\Gamma_{W}(T,\Z)\simeq R\Gamma_{Wh}(T,\Z)$. 
\end{prop}

\begin{proof}
Since \'etale cohomology with coefficients in $\Z$ does not change
if we replace $T$ by $T^\red$, we can assume that $T$ is reduced.
We proceed by induction on dimension of $T$ and the number of irreducible 
components of $T$.
If the number of components is one, then $T^\red$ is smooth and proper
and the result follows from \cite[Thm. 4.3]{Geisser06}.
In general, let $T=\cup_{i\in I}S_i$ and set $Y=S_1$ and $T'=\cup_{i\not=1}S_i$.
Then the hypothesis of Lemma \ref{finiteblowup} are satisfied, $Y$ is smooth,
$T'$ is a normal crossing scheme with fewer irreducible components, and
$Y'$ a normal crossing scheme of smaller dimension.
Hence we obtain the first statement on perfectness, and the second statement by
comparing with the corresponding triangle for Wh-cohomology.
\end{proof}

Note that we can have $H_{W}^{2}(T,\Z)\not= H_{Wh}^{2}(T,\Z)$ for normal
proper surfaces \cite[Prop. 8.2]{Geisser06}.

\subsection{Finite generation of homology}
For later use, we record the two following conditional results on 
finite generation.

To get statements for homology, recall the following conjecture 
from \cite{Geisser10b}.

\medskip
\noindent{\bf Conjecture} 
$P_n(X)$: For the smooth and proper scheme $X$ over a finite field, 
the group $CH_n(X, i)$ is torsion for all $i > 0$.
\medskip

Conjecture $P_n(X)$ is known for all $n$ if $X$ is a curve.
In general, it is a particular case of Parshin's conjecture, which is equivalent
to the statement $P_n(X)$ for all $n$.  
Parshin's conjecture in turn is implied by the Beilinson-Tate conjecture
\cite[Thm.\ 1.2]{Geisser98}. By the projective bundle formula,
conjecture $P_n(X)$ for all $X$ of dimension at most $d$ 
implies conjecture $P_{n-1}(X)$ for all $X$ of dimension $d-1$. 
The following Proposition is 
\cite[Prop. 4.2]{Geisser10b}.

\begin{prop}\label{homology.fin.gen}
If conjecture $P_{0}(X)$ holds for all smooth and proper 
schemes of dimension at
most $\dim Y$, then the cohomology groups of 
$R\Gamma_W(Y,\Z^c(0))$ are finitely generated
and vanish for almost all indices. 
\end{prop}

If $Y$ is a simple normal crossing scheme, then it suffices to assume 
that $P_{0}(Y_I)$ holds for all multiple intersections $Y_I$.


\begin{prop}\label{-1condition}
If resolution of singularities and conjecture $P_{-1}(X)$ holds for 
all schemes of dimension at most $d-1$, then the cohomology groups of 
$R\Gamma_{Wh}(\mathcal{X}_s,\Z(d))$
are finite and vanish for almost all indices. 
\end{prop}

\begin{proof}
Using blow-up squares and induction on the dimension it suffices to
prove the statement for smooth and proper schemes $T$  of dimension 
at most $d-1$. By \cite[Cor. 5.5]{Geisser06} the Weil-eh cohomology
groups agree with Weil-\'etale cohomology groups. By conjecture 
$P_{-1}(T)$ they are torsion, hence finite by comparison with 
\'etale cohomology groups.
\end{proof}

\subsection{Duality}\label{sectdualityXs}

\begin{thm}\label{dualityXs}
Let $Y$ be a simple normal crossing scheme of 
over a finite field $k$ such that $R\Gamma_W(Y,\Z^c(0))$ is a perfect complex of abelian groups. Then there is  a perfect pairing
\begin{equation}\label{pairingperf}
R\Gamma_{W}(Y,\Z)\otimes^\bL R\Gamma_W(Y,\Z^c(0))\rightarrow  \Z[-1]
\end{equation}
of perfect complexes of abelian groups.
\end{thm}

\begin{proof}
Let $f:Y\rightarrow\mathrm{Spec}(k)=s$ be the structure morphism. 
The push-forward map \cite[Cor.\ 3.2]{Geisser10}
$$Rf_*\Z^c(0)^{Y}\rightarrow \Z^c(0)^{s}\simeq\Z[0]$$
induces a trace map
$$R\Gamma_{W}(Y,\Z^c(0))\rightarrow R\Gamma_{W}(s,\Z)
\rightarrow \Z[-1].$$
We consider the map
\begin{equation}\label{pairingperfind}
R\Gamma_{W}(Y,\Z)\longrightarrow R\mathrm{Hom}(R\Gamma_W(Y,\Z^c(0)), \Z[-1])
\end{equation}
induced by the pairing
\begin{equation*}\label{pairingperf*}
R\Gamma_{W}(Y,\Z)\otimes^{\bL}R\Gamma_W(Y,\Z^c(0))\rightarrow 
R\Gamma_W(Y,\Z^c(0)) \rightarrow R\Gamma_W(s,\Z^c(0))\rightarrow \Z[-1]
\end{equation*}
which in turn is induced by the obvious pairing  
$\Z\otimes^{\bL}\Z^c(0)\rightarrow \Z^c(0)$.
In order to show that the morphism of perfect complexes (\ref{pairingperfind}) 
is an equivalence, it is enough to show that 
$(\ref{pairingperfind})\otimes^{\bL}\Z/m\Z$ is an equivalence for any 
integer $m$. But $(\ref{pairingperfind})\otimes^{\bL}\Z/m\Z$ may be 
identified with the canonical map
\begin{equation}\label{pairinghere}
R\Gamma_{et}(Y,\Z/m\Z)\longrightarrow 
R\mathrm{Hom}(R\Gamma_{et}(Y,\Z^c(0)/m), \bq/\Z[-1])
\end{equation}
since we have an equivalence of lax symmetric monoidal functors 
$$R\Gamma_W(Y,-)\otimes^{\bL}_{\Z}\Z/m\Z\simeq 
R\Gamma_{et}(Y,(-)\otimes^{\bL}_{\Z}\Z/m\Z).$$
But (\ref{pairinghere}) is an equivalence by \cite[Theorem 5.1]{Geisser10}. 
Hence (\ref{pairingperfind}) is an equivalence as well. 
\end{proof}

\section{The complexes $R\Gamma_{ar}(\mathcal{X}_K,\Z(n))$ 
in $\mathbf{D}^b(\mathrm{LCA})$.}
Under the assumption that the pair $(\mathcal{X},n)$ satisfies Hypothesis 
\ref{hyp-reductionmap} below,  we give in Section \ref{sectgeneral} a 
construction of complexes in a fiber sequence in $\mathbf{D}^b(\mathrm{LCA})$
\begin{equation}\label{introsecttriangle}
R\Gamma_{ar}(\mathcal{X}_s,Ri^{!}\Z(n))\rightarrow R\Gamma_{ar}(\mathcal{X},\Z(n))\rightarrow R\Gamma_{ar}(\mathcal{X}_K,\Z(n))
\end{equation}
where 
$$R\Gamma_{ar}(\mathcal{X}_s,Ri^{!}\Z(n)):=R\Gamma_{W}(\mathcal{X}_s,Ri^{!}\Z(n))$$
is defined below. Hypothesis \ref{hyp-reductionmap} is known for $n=0,1$ and arbitrary $\mathcal{X}$, hence this construction is unconditional in those cases. 

In Sections \ref{uncondsect0} and \ref{uncondsectd} we give an alternative definition of the triangle \eqref{introsecttriangle} for $n=0$ and 
$n=d:=\mathrm{dim}(\mathcal{X})$, respectively, which are expected to coincide 
with the conditional definition of Section \ref{sectgeneral}. 

In Section \ref{sectionfiniteranks} we show that these complexes in fact belong to $\mathbf{D}^b(\mathrm{FLCA})$ under some conditions. In Section \ref{LCAcohgroups}, we show that the cohomology of these complexes consists of locally compact abelian groups for $n=0,d$.

\subsection{Notation}\label{sectnotgeneral}
Let $p$ be a prime number, let $K/\mathbb{Q}_p$ be a  finite extension, let $\mathcal{O}_K$ be its ring of integers and let $\bar{K}/K$ be an algebraic closure. We denote by $K^{un}$ the maximal unramified extension of $K$ inside $\bar{K}$. Let $\mathcal{X}/\mathcal{O}_K$ be a regular, proper and flat scheme over $\mathrm{Spec}(\mathcal{O}_K)$. Suppose that $\X$ is connected of Krull dimension $d$. Let $\mathcal{X}_s$ be its special fiber, where $s\in\mathrm{Spec}(\mathcal{O}_K)$ is the closed point. We consider the following diagram.
\[ \xymatrix{
\mathcal{X}_{K^{un}}\ar[r]^{\bar{j}}\ar[d]^{}&\mathcal{X}_{\mathcal{O}_{K^{un}}}\ar[d]& \ar[l]_{\bar{i}}\mathcal{X}_{\bar{s}}\ar[d]\\
\mathcal{X}_{K} \ar[r]^{j} &\mathcal{X}&\ar[l]_{i}\mathcal{X}_{s} .
}
\]
For any $n\geq0$, we denote by $\Z(n)$ Bloch's cycle complex considered 
as a complex of \'etale sheaves and 
$\Z/m(n):=\Z(n)\otimes^\bL\Z/m$
We denote  by $G_{\kappa(s)}\simeq\widehat{\Z}$ and by $W_{\kappa(s)}\simeq\Z$ the Galois group and the Weil group of the finite field $\kappa(s)$, respectively. We define Weil-\'etale cohomology groups
\begin{eqnarray*}
R\Gamma_W(\mathcal{X}_K,\Z(n))&:=&R\Gamma(W_{\kappa(s)},R\Gamma_{et}(\mathcal{X}_{K^{un}},\Z(n))),\\
R\Gamma_W(\mathcal{X},\Z(n))&:=&R\Gamma(W_{\kappa(s)},R\Gamma_{et}(\mathcal{X}_{\mathcal{O}_{K^{un}}},\Z(n))),\\
R\Gamma_{W}(\mathcal{X}_{s},Ri^!\Z(n))&:=&R\Gamma(W_{\kappa(s)},R\Gamma_{et}(\mathcal{X}_{\bar{s}},R\bar{i}^!\Z(n)).
\end{eqnarray*}
Applying $R\Gamma(W_{\kappa(s)},-)$ to 
the fiber sequence
\begin{equation*}
R\Gamma_{et}(\mathcal{X}_{\bar{s}},R\bar{i}^!\Z(n))\rightarrow R\Gamma_{et}(\mathcal{X}_{\mathcal{O}_{K^{un}}},\Z(n))\rightarrow R\Gamma_{et}(\mathcal{X}_{K^{un}},\Z(n)),
\end{equation*}
we obtain the fiber sequence
\begin{equation}\label{locfibsequenceW}
R\Gamma_{W}(\mathcal{X}_{s},Ri^!\Z(n))\rightarrow R\Gamma_W(\mathcal{X},\Z(n))\rightarrow R\Gamma_W(\mathcal{X}_K,\Z(n)).
\end{equation}

\subsection{Uniform conditional definition}\label{sectgeneral}
Recall from \cite{Geisser06} the $eh$-motivic cohomology $R\Gamma_{eh}(-,\Z(n))$
and $Wh$-motivic cohomology
$$R\Gamma_{Wh}(\mathcal{X}_{s},\Z(n)):=R\Gamma(W_{\kappa(s)},R\Gamma_{eh}(\mathcal{X}_{\overline{s}},\Z(n))).$$

\begin{Hypothesis}\label{hyp-reductionmap}
We have a reduction map 
$$\bar{i}^{*}:R\Gamma_{et}(\mathcal{X}_{\mathcal{O}_{K^{un}}},\Z(n))\rightarrow R\Gamma_{eh}(\mathcal{X}_{\bar{s}},\Z(n))$$
and the complexes $R\Gamma_{et}(\mathcal{X},\Z(n))$, 
$R\Gamma_{eh}(\mathcal{X}_{s},\Z(n))$ and 
$R\Gamma_{et}(\mathcal{X}_{s},Ri^!\Z(n))$ are cohomologically bounded.
\end{Hypothesis}

\begin{defn}
Under hypothesis \ref{hyp-reductionmap}, we apply the functor $R\Gamma(W_{\kappa(s)},-)$ to the reduction map $\bar{i}^{*}$ and we obtain a map
\begin{equation}\label{mapW}
R\Gamma_{W}(\mathcal{X},\Z(n))\rightarrow R\Gamma_{Wh}(\mathcal{X}_{s},\Z(n)).
\end{equation}
We denote the cofiber of (\ref{mapW}) by $C_W(\mathcal{X},n)$, 
so that we have a cofiber sequence
\begin{equation}\label{cofibW}
R\Gamma_{W}(\mathcal{X},\Z(n))\rightarrow R\Gamma_{Wh}(\mathcal{X}_{s},\Z(n))\rightarrow C_W(\mathcal{X},n)
\end{equation}
 in $\mathbf{D}^b(\mathrm{Ab})$.
\end{defn}

\begin{prop}\label{topologicalar}
Assume Hypothesis \ref{hyp-reductionmap}. Then there exist $\mathbf{R\Gamma}_{ar}(\mathcal{X},\Z(n))\in\mathbf{D}^b(\mathrm{LCA})$ and $\mathbf{R\Gamma}_{ar}(\mathcal{X}_K,\Z(n))\in\mathbf{D}^b(\mathrm{LCA})$ endowed with fiber sequences
\begin{equation}\label{fibseqar1}
\mathbf{R\Gamma}_{ar}(\mathcal{X},\Z(n))\rightarrow R\Gamma_{Wh}(\mathcal{X}_s,\Z(n))\rightarrow C_W(\mathcal{X},n)\widehat{\underline{\otimes}}\widehat{\Z}
\end{equation}
and
\begin{equation}\label{fibseqar2}
R\Gamma_{W}(\mathcal{X}_s,Ri^{!}\Z(n))\rightarrow \mathbf{R\Gamma}_{ar}(\mathcal{X},\Z(n))\rightarrow \mathbf{R\Gamma}_{ar}(\mathcal{X}_K,\Z(n))
\end{equation}
in $\mathbf{D}^b(\mathrm{LCA})$.
\end{prop}
\begin{proof}
Composing the morphism in $\mathbf{D}^b(\mathrm{Ab})$
$$R\Gamma_{Wh}(\mathcal{X}_s,\Z(n))\rightarrow C_W(\mathcal{X},n)$$
and
the morphism in $\mathbf{D}^b(\mathrm{LCA})$
$$C_W(\mathcal{X},n)\rightarrow C_W(\mathcal{X},n)\widehat{\underline{\otimes}}\widehat{\Z}$$
given by Lemma \ref{lemmaptohat}, we obtain a morphism in $\mathbf{D}^b(\mathrm{LCA})$
\begin{equation}\label{maphere}
R\Gamma_{Wh}(\mathcal{X}_s,\Z(n))\rightarrow C_W(\mathcal{X},n)\widehat{\underline{\otimes}}\widehat{\Z}.
\end{equation}
We define $\mathbf{R\Gamma}_{ar}(\mathcal{X},\Z(n))$ as the fiber of (\ref{maphere}) and we obtain the fiber sequence (\ref{fibseqar1}) in $\mathbf{D}^b(\mathrm{LCA})$.
Lemma \ref{lemmaptohat} gives a  map from (\ref{cofibW}) to (\ref{fibseqar1}) hence a map
$$R\Gamma_{W}(\mathcal{X},\Z(n))\rightarrow \mathbf{R\Gamma}_{ar}(\mathcal{X},\Z(n)).$$
Then we define $\mathbf{R\Gamma}_{ar}(\mathcal{X}_K,\Z(n))\in \mathbf{D}^b(\mathrm{LCA})$ as the cofiber of the composite map
$$R\Gamma_{W}(\mathcal{X}_s,Ri^{!}\Z(n))\rightarrow R\Gamma_{W}(\mathcal{X},\Z(n))\rightarrow \mathbf{R\Gamma}_{ar}(\mathcal{X},\Z(n)).$$
\end{proof}

\begin{rem}
Since $\Z(0)\cong \Z$ and $\Z(1)\cong {\mathbb G}_m[-1]$, 
Hypothesis \ref{hyp-reductionmap} holds for $n=0$ and $n=1$, 
so that $\mathbf{R\Gamma}_{ar}(\mathcal{X},\Z(n))$ and $\mathbf{R\Gamma}_{ar}(\mathcal{X}_K,\Z(n))$ are unconditionally defined in these cases.
\end{rem}

\subsection{Working definition for the Tate twist $n=0$.}\label{uncondsect0}

We assume that $\X^{\red}_s$ is a simple normal crossing scheme.
To obtain unconditional definitions for $n=0$, we replace 
$R\Gamma_{Wh}(\mathcal{X}_s,\Z)$ by $R\Gamma_{W}(\mathcal{X}_s,\Z)$
in the construction of Section \ref{sectgeneral}. 
In view of Corollary \ref{corW*}, this will agree with the definition of 
Section \ref{sectgeneral} provided that resolution of singularities for 
schemes of dimension at most $\mathrm{dim}(\X_s)$ exist.

There is a canonical map
$$R\Gamma_{W}(\mathcal{X},\Z)\rightarrow 
R\Gamma_{W}(\mathcal{X}_s,\Z)$$
whose cofiber we again denote by $C_W(\mathcal{X},0)$. Following the construction of Section \ref{sectgeneral}, we define $R\Gamma_{ar}(\mathcal{X},\Z)\in\mathbf{D}^b(\mathrm{LCA})$ and $R\Gamma_{ar}(\mathcal{X}_K,\Z)\in\mathbf{D}^b(\mathrm{LCA})$ endowed with fiber sequences
\begin{equation}\label{fibseqar1*}
R\Gamma_{ar}(\mathcal{X},\Z)\rightarrow R\Gamma_{W}(\mathcal{X}_s,\Z)\rightarrow C_W(\mathcal{X},0)\widehat{\underline{\otimes}}\widehat{\Z}
\end{equation}
and
\begin{equation}\label{fibseqar2*}
R\Gamma_{W}(\mathcal{X}_s,Ri^{!}\Z)\rightarrow R\Gamma_{ar}(\mathcal{X},\Z)\rightarrow R\Gamma_{ar}(\mathcal{X}_K,\Z)
\end{equation}
in $\mathbf{D}^b(\mathrm{LCA})$. We used bold letters for the complexes defined in Section \ref{sectgeneral} in order to distinguish them with the complexes defined in this Section.

\begin{prop}\label{red-map-n=0}
If $\X^{\red}_s$ is a simple normal crossing  scheme, then the map
$$R\Gamma_{ar}(\mathcal{X},\Z)\rightarrow 
R\Gamma_{W}(\mathcal{X}_s,\Z)$$
is an  equivalence. For arbitrary $\X$, the map
$$\mathbf{R\Gamma}_{ar}(\mathcal{X},\Z)\rightarrow 
R\Gamma_{Wh}(\mathcal{X}_s,\Z)$$
is an  equivalence.
\end{prop}

\begin{proof}
By proper base change, the map
$$R\Gamma_{et}(\mathcal{X},\Z/m)\rightarrow R\Gamma_{et}(\mathcal{X}_s,\Z/m)$$
is an equivalence, hence $C_W(\mathcal{X},0)\otimes^{\bL}\Z/m\simeq 0$.
We obtain $C_W(\mathcal{X},0)\underline{\widehat{\otimes}}\widehat{\Z}\simeq 0$. The first equivalence of the proposition follows. The second equivalence is obtained the same way, in view of the fact that
$$R\Gamma_{et}(\mathcal{X},\Z/m)\rightarrow R\Gamma_{et}(\mathcal{X}_s,\Z/m)\rightarrow R\Gamma_{eh}(\mathcal{X}_s,\Z/m)$$
is an equivalence, again by proper base change.
\end{proof}

\begin{prop}\label{samedefn=0}
If $\X^{\red}_s$ is a simple normal crossing  scheme, then there is a canonical map of fiber sequences
\[ \xymatrix{
R\Gamma_{W}(\mathcal{X}_s,Ri^{!}\Z)\ar[r]^{}\ar@{=}[d]&R\Gamma_{ar}(\mathcal{X},\Z)\ar[d] \ar[r]& R\Gamma_{ar}(\mathcal{X}_K,\Z)\ar[d]\\
R\Gamma_{W}(\mathcal{X}_s,Ri^{!}\Z)\ar[r]&\mathbf{R\Gamma}_{ar}(\mathcal{X},\Z) \ar[r]& \mathbf{R\Gamma}_{ar}(\mathcal{X}_K,\Z)
}
\]
If resolution of singularities for schemes over $\kappa(s)$ of dimension at most 
$ d-1$ \cite[Definition 2.4]{Geisser06} exist, then this morphism of fiber sequences is an equivalence. 
\end{prop}

\begin{proof}
This follows from Corollary \ref{corW*} and Proposition \ref{red-map-n=0}.
\end{proof}

In particular, we obtain that 
$$R\Gamma_{ar}(\mathcal{X},\Z)\stackrel{\sim}{\longrightarrow} 
\mathbf{R\Gamma}_{ar}(\mathcal{X},\Z)$$
is an equivalence if $d\leq 3$.

\begin{notation}\label{notZcoeff}
If $\X^{\red}_s$ is a simple normal crossing  scheme, we denote by 
$R\Gamma_{ar}(\X,\Z)$ and $R\Gamma_{ar}(\X_K,\Z)$ the complexes defined above. 
In view of Proposition \ref{samedefn=0}, we set 
$R\Gamma_{ar}(\X,\Z):=\mathbf{R\Gamma}_{ar}(\X,\Z)$ and 
$R\Gamma_{ar}(\X_K,\Z):=\mathbf{R\Gamma}_{ar}(\X_K,\Z)$ 
for arbitrary regular $\X$ of dimension at most $3$ or when we are
assuming resolution of singularities.
\end{notation}

\subsection{Working definition for the Tate twist $n=d$.}\label{uncondsectd}

The complex $R\Gamma_W(\mathcal{X},\Z(d))$ is not known to be bounded below. 
However the complex $$R\Gamma_W(\mathcal{X},\mathbb{Q}/\Z(d))
\simeq R\Gamma_{et}(\mathcal{X},\mathbb{Q}/\Z(d))$$
is bounded, as can be seen by duality, hence the cohomology groups 
$H^{i}_W(\mathcal{X},\Z(d))$ are $\bq$-vector spaces for $i\ll 0$. In particular, for $a<b\ll 0$ the map
$$\tau^{>a}R\Gamma_W(\mathcal{X},\Z(d))\rightarrow \tau^{>b}R\Gamma_W(\mathcal{X},\Z(d))$$
induces an equivalence
$$(\tau^{>a}R\Gamma_W(\mathcal{X},\Z(d)))\widehat{\underline{\otimes}}\widehat{\Z}\stackrel{\sim}{\longrightarrow} (\tau^{>b}R\Gamma_W(\mathcal{X},\Z(d)))\widehat{\underline{\otimes}}\widehat{\Z}.$$

\begin{defn}\label{uncdefn=d}
Let $a\ll 0$. We define $$R\Gamma_{ar}(\mathcal{X},\Z(d)):= (\tau^{>a}R\Gamma_{W}(\mathcal{X},\Z(d)))\widehat{\underline{\otimes}}\widehat{\Z}.$$
If $R\Gamma_W(\mathcal{X}_s,\Z^c(0))$ is cohomologically bounded, we define $R\Gamma_{ar}(\mathcal{X}_K,\Z(d))$ as the cofiber of the composite map
$$R\Gamma_W(\mathcal{X}_s,Ri^!\Z(d))\rightarrow \tau^{>a}R\Gamma_{W}(\mathcal{X},\Z(d))\rightarrow R\Gamma_{ar}(\mathcal{X},\Z(d))$$
in $\mathbf{D}^b(\mathrm{LCA})$. 
\end{defn}

\begin{rem}\label{rempurityshifted}
On the connected, $d$-dimensional and regular scheme $\X$, we have  $\Z(d)^{\X}=\Z^c(0)^{\X}[-2d]$ by definition. By \cite[Corollary 7.2]{Geisser10}, we have $Ri^!\Z^c(0)^{\X}=\Z^c(0)^{\X_s}$, hence $Ri^!\Z(d)^{\X}=\Z^c(0)^{\X_s}[-2d]$. 
\end{rem}

\begin{prop}\label{propTB}
Suppose that $\X$ satisfies Hypothesis \ref{hyp-reductionmap} for $n=d$
and suppose that $R\Gamma_{W}(\mathcal{X}_K,\Z(d))$ and
$R\Gamma_W(\mathcal{X}_s,\Z^c(0))$ are cohomologically bounded. Then there is a canonical map of fiber sequences
\[ \xymatrix{
R\Gamma_{W}(\mathcal{X}_s,Ri^{!}\Z(d)) \ar[r]^{}\ar@{=}[d]&\mathbf{R\Gamma}_{ar}(\mathcal{X},\Z(d))\ar[d] \ar[r]& \mathbf{R\Gamma}_{ar}(\mathcal{X}_K,\Z(d))\ar[d]\\
R\Gamma_{W}(\mathcal{X}_s,Ri^{!}\Z(d))\ar[r]&R\Gamma_{ar}(\mathcal{X},\Z(d)) \ar[r]& R\Gamma_{ar}(\mathcal{X}_K,\Z(d)).
}
\]
If $R\Gamma_{Wh}(\mathcal{X}_s,\Z(d))$ has finite cohomology groups, then this 
morphism of fiber sequences is an equivalence.
\end{prop}

Cohomological boundedness of  
$R\Gamma_{W}(\mathcal{X}_K,\Z(d))$ in negative degrees 
is a special case of 
the Beilinson-Soul\'e conjecture stating that there is no
negative motivic cohomology, and in positive degrees it follows
for finite cohomological dimension reasons. 
See Proposition \ref{homology.fin.gen} and
Proposition \ref{-1condition} for the other boundedness conditions.

\begin{proof}
If $R\Gamma_{W}(\mathcal{X}_K,\Z(d))$ and
$R\Gamma_{W}(\mathcal{X}_s,Ri^!\Z(d))\simeq R\Gamma_W(\mathcal{X}_s,\Z^c(0))[-2d]$ are cohomologically bounded,
then the same holds for $R\Gamma_{W}(\mathcal{X},\Z(d))$
by the localization triangle
$$ \cdots\to R\Gamma_{W}(\mathcal{X}_s,Ri^!\Z(d))
\to R\Gamma_{W}(\mathcal{X},\Z(d))\to 
R\Gamma_{W}(\mathcal{X}_K,\Z(d))\to \cdots.$$
In this case, 
$R\Gamma_{W}(\mathcal{X},\Z(d))\stackrel{\sim}{\rightarrow} \tau^{>a}R\Gamma_{W}(\mathcal{X},\Z(d))$ is an equivalence for $a<<0$, hence
$$R\Gamma_{ar}(\mathcal{X},\Z(d))\simeq R\Gamma_{W}(\mathcal{X},\Z(d))\widehat{\underline{\otimes}}\widehat{\Z}.$$
In view of (\ref{cofibW}) and (\ref{fibseqar1}), we obtain a map of fiber sequences
\[ \xymatrix{
\mathbf{R\Gamma}_{ar}(\mathcal{X},\Z(d))\ar[r]^{}\ar[d]^{}&R\Gamma_{Wh}(\mathcal{X}_s,\Z(d))\ar[d] \ar[r]& C_W(\X,d)\widehat{\underline{\otimes}}\widehat{\Z}\ar@{=}[d]\\
R\Gamma_{ar}(\mathcal{X},\Z(d))\ar[r]&R\Gamma_{Wh}(\mathcal{X}_s,\Z(d)) \widehat{\underline{\otimes}}\widehat{\Z}\ar[r]& C_W(\X,d)\widehat{\underline{\otimes}}\widehat{\Z}
}
\]
since $(-)\widehat{\underline{\otimes}}\widehat{\Z}$ is an exact functor. 
If the cohomology of $R\Gamma_{Wh}(\mathcal{X}_s,\Z(d))$ consists of finite groups,
then the middle vertical map 
$$R\Gamma_{Wh}(\mathcal{X}_s,\Z(d))\rightarrow R\Gamma_{Wh}(\mathcal{X}_s,\Z(d)) \widehat{\underline{\otimes}}\widehat{\Z}$$
in the above diagram is an equivalence by Proposition \ref{discprocompletion},
and the two fiber sequences are equivalent.
\end{proof}

\begin{rem}\label{rem-ar-model=et}
It follows from Proposition \ref{discprocompletion} and Remark \ref{remfordiscrete} that we have 
\begin{eqnarray*}
\mathrm{disc} (R\Gamma_{ar}(\X,\Z(d)))&\simeq & R\mathrm{lim}\,(\tau^{>a}R\Gamma_{W}(\mathcal{X},\Z(d))\otimes^{\bL}\Z/m)\\
&\simeq& R\mathrm{lim}\,R\Gamma_{W}(\mathcal{X},\Z/m(d))\\
&\simeq&R\mathrm{lim}\,R\Gamma_{et}(\X,\Z/m(d))\\
&=:& R\Gamma_{et}(\X,\widehat{\Z}(d)) 
\end{eqnarray*}
where we denote $\Z/m(d):=\Z(d)\otimes^{\bL}\Z/m$.
\end{rem}

\subsection{Finite ranks}\label{sectionfiniteranks}

\begin{lem}\label{finranklemma}
If $\mathcal{X}$ is normal, then we have an isomorphism 
$$H^j_{et}(\mathcal{X}_s,Ri^!\mathbb{Q}/\Z)\cong 
H^{j+1}_{et}(\mathcal{X}_s,Ri^!\Z)$$
of abelian groups of finite ranks for all $j\in\Z$. 
\end{lem}
\begin{proof}
The isomorphism follows because $\X$ normal implies that 
$\mathbb{Q}\cong Rj_*j^*\mathbb{Q}$, hence $Ri^!\mathbb{Q}\cong 0$.
Since $H^j_{et}(\mathcal{X}_s,Ri^!\mathbb{Q}/\Z)$ is torsion and discrete, 
it is both of finite $\Z$-rank and of finite $\mathbb{S}^1$-rank. 
It remains to see that it is of finite $p$-rank for any prime number $p$.
But $H^j_{et}(\mathcal{X}_s,Ri^!\Z/p\Z)$ is a finite group for any $j\in\Z$, because of the fiber sequence
$$R\Gamma_{et}(\mathcal{X}_s,Ri^!\Z/p\Z)\rightarrow 
R\Gamma_{et}(\mathcal{X},\Z/p\Z) \rightarrow 
R\Gamma_{et}(\mathcal{X}_K,\Z/p\Z)$$
and classical finiteness results in  \'etale and Galois cohomology.
\end{proof}

\begin{prop} \label{prop-finiteranks}
a) Assume that $\X_s^{\red}$ is a simple normal crossing  scheme or
assume resolution of singularities for schemes over $\kappa(s)$ of dimension at
most $d-1$ \cite[Definition 2.4]{Geisser06}. Then 
$R\Gamma_{ar}(\mathcal{X},\Z)$ and 
$R\Gamma_{ar}(\mathcal{X}_K,\Z)$ belong to $\mathbf{D}^b(\mathrm{FLCA})$.

b)  Assume that 
$R\Gamma_W(\mathcal{X}_s,\Z^c(0))$ is a perfect complex of abelian groups. 
Then $R\Gamma_{ar}(\mathcal{X},\Z(d))$ and 
$R\Gamma_{ar}(\mathcal{X}_K,\Z(d))$ belong to $\mathbf{D}^b(\mathrm{FLCA})$.
\end{prop}

\begin{proof}
a) Under the hypothesis, the complexes $R\Gamma_{W}(\mathcal{X}_s,\Z)$
and $R\Gamma_{Wh}(\mathcal{X}_s,\Z)$ 
are perfect complexes of abelian groups by Proposition \ref{prop-normalcrossing} 
and Proposition \ref{propperfectWh}, respectively, 
hence they belong to $\mathbf{D}^b(\mathrm{FLCA})$ by Lemma \ref{lemcoh}. 
The result for $R\Gamma_{ar}(\mathcal{X},\Z)$
then follows from Proposition \ref{red-map-n=0} (using Notation \ref{notZcoeff}),
and the result for $R\Gamma_{ar}(\mathcal{X}_K,\Z)$ follows from
Proposition \ref{samedefn=0} and Lemma \ref{finranklemma}.

b) By the proof of Proposition \ref{prop-localduality}, 
$R\Gamma_{ar}(\mathcal{X},\Z(d))$ is (up to a shift) dual to 
$R\Gamma_{et}(\mathcal{X}_s,Ri^!\mathbb{Q}/\Z)$. 
Hence the result follows from Lemmas \ref{lemcoh} and \ref{finranklemma}. 
The statement for $R\Gamma_{ar}(\mathcal{X}_K,\Z(d))$ follows from
the statement $R\Gamma_{ar}(\mathcal{X},\Z(d))$ together with the 
perfectness of 
$R\Gamma_W(\mathcal{X}_s,Ri^!\Z(d))\simeq R\Gamma_W(\mathcal{X}_s,\Z^c(0))[-2d]$ by hypothesis.
\end{proof}

\subsection{The topology on cohomology groups}\label{LCAcohgroups}

Recall from Section \ref{sectnotgeneral} that $\X$ denotes a regular connected scheme which is proper and flat over $\mathcal{O}_K$.  We refer to \cite[Section 2.16]{deJong96} for the following definition.

\begin{defn}
Let $\X_{s,i}, i\in I$ be the irreducible components of $\X_s$. We set $\X_{s,J}=\cap_{i\in J} \X_{s,i}$ for any non-empty subset $J\subseteq I$. We say that $\X/\mathcal{O}_K$ has \emph{strictly semi-stable reduction} if $\X_s$ is reduced, $\X_{s,i}$ is a divisor on $\X$, and for each non-empty $J\subseteq I$, the scheme $\X_{s,J}$ is smooth over $\kappa(s)$ and has codimension $\vert J\vert$ in $\X$.

\end{defn}
If $\X/\mathcal{O}_K$ has strictly semi-stable reduction then $\X_s$ is a simple normal crossing scheme over $\kappa(s)$, in the sense of Definition \ref{strictNCS}.

\begin{thm}\label{leminj}
Suppose that $\X/\mathcal{O}_K$ has strictly semi-stable reduction. Then for any $i\in\Z$, the map
\begin{equation}\label{mapSato}
H^{i}_{et}(\mathcal{X},\mathbb{Q}_p(d))\rightarrow H^{i}_{et}(\mathcal{X}_K,\mathbb{Q}_p(d))
\end{equation}
 is injective.
\end{thm}

\begin{proof}
Since $\X/\mathcal{O}_K$ has strictly semi-stable reduction, the morphism 
$\X\rightarrow\mathrm{Spec}(\mathcal{O}_K)$ is log smooth with respect to the 
log structures associated with $\X_s$ and $s$ respectively, 
where $s$ is the closed point of $\mathrm{Spec}(\mathcal{O}_K)$, and 
$\X_s$ is a normal crossing divisor on $\X$. 
Therefore, the results of \cite{Sato20} apply. We have isomorphisms
\begin{equation}\label{isoSato}
H^{i}_{et}(\mathcal{X},\mathbb{Q}_p(d))\simeq H^{i}_{et}(\mathcal{X},\mathbb{Q}^S_p(d))
\end{equation}
compatible with the map (\ref{mapSato}), where 
$$R\Gamma_{et}(\mathcal{X},\mathbb{Q}^S_p(d)):=
R\mathrm{lim}R\Gamma_{et}(\mathcal{X},\mathfrak{T}_r(d))\otimes_{\Z_p}\bq_p$$ 
is the complex studied in \cite{Sato07} and \cite{Sato20}. 
Indeed, this follows from the equivalences 
\begin{eqnarray}
R\Gamma(\X_{et},\mathfrak{T}_r(d))&\simeq& R\mathrm{Hom}(R\Gamma_{\X_s}(\X_{et},\Z/p^{r}), \Z/p^r)[-2d-1]\\
&\simeq&  R\Gamma(\X_{et},\Z/p^r(d)),
\end{eqnarray}
given by \cite[Thm 10.1.1]{Sato07} and \cite[Proof of Thm 7.5]{Geisser10}, and from the fact that (\ref{mapSato}) is induced by the dual of the map
$$R\Gamma(\X_{K,et},\Z/p^{r})[-1]\rightarrow R\Gamma_{\X_s}(\X_{et},\Z/p^{r}).$$
Hence we are reduced to show that the map
$$H^{i}_{et}(\mathcal{X},\mathbb{Q}^S_p(d))\rightarrow H^{i}_{et}(\mathcal{X}_K,\mathbb{Q}_p(d))$$
is injective. By \cite[Prop. 3.4(1)]{Sato20}, \cite[Section 4.1]{Sato20} and \cite[Thm. 5.3]{Sato20}, there is a morphism of spectral sequences from
$$H^{i}_f(G_K, H^{j}_{et}(\X_{\bar{K}},\bq_p(d)))\Rightarrow H^{i+j}_{et}(\mathcal{X},\mathbb{Q}^S_p(d))$$ 
to 
$$H^{i}(G_K, H^{j}_{et}(\X_{\bar{K}},\bq_p(d)))\Rightarrow H^{i+j}_{et}(\mathcal{X}_K,\mathbb{Q}_p(d))$$
where the first spectral sequence degenerates into isomorphisms
$$H^{j}_{et}(\mathcal{X},\mathbb{Q}^S_p(d))\stackrel{\sim}{\rightarrow}H^{1}_f(G_K, H^{j-1}_{et}(\X_{\bar{K}},\bq_p(d))).$$
Since we have \cite[Prop. 5.10(1)]{Sato20} 
$$H^{0}_f(G_K, H^{j}_{et}(\X_{\bar{K}},\bq_p(d)))=H^{0}(G_K, H^{j}_{et}(\X_{\bar{K}},\bq_p(d)))=0$$
for any $j\in\Z$, we obtain a commutative square
\[ \xymatrix{
H^{j}_{et}(\mathcal{X},\mathbb{Q}^S_p(d))  \ar[r]\ar[d]^{\simeq} &  H^{j}_{et}(\mathcal{X}_K,\mathbb{Q}_p(d))\ar[d]\\
H^{1}_f(G_K, H^{j-1}_{et}(\X_{\bar{K}},\bq(d)))  \ar[r]& H^{1}(G_K, H^{j-1}_{et}(\X_{\bar{K}},\bq_p(d)))
}
\]
where the vertical maps are edge morphisms of the corresponding spectral sequences. Here the left vertical map is an isomorphism and the lower horizontal map is injective. It follows that the upper horizontal map is injective as well.
\end{proof}

\begin{lem}\label{lemfiniteness}
Suppose that $\X_s$ is a simple normal crossing scheme. Then the group
$$H^{i}_{et}(\mathcal{X},\widehat{\Z}'(d)):= 
H^{i}(R\mathrm{lim}_{p\nmid m} \,
R\Gamma_{et}(\mathcal{X},\Z(d))\otimes^{\bL}\Z/m)$$
is finite for any $i\in\Z$.
\end{lem}

\begin{proof}
By definition we have $\Z(d)\cong \Z^c(0)[-2d]$, and by 
\cite[Prop. 7.10 a)]{Geisser10} and Gabber's purity theorem
\cite[\S 8]{Fujiwara} we have 
$$\Z^c/m(0)[-2d] \cong Rf^!\Z^c/m(0)[-2d]\cong 
Rf^!\Z/m[-2d]\cong \mu_m^{\otimes d}$$
on $\X$ for any $m$ prime to $p$. Moreover, by the proper base change theorem
$$R\Gamma_{et}(\X,\mu_m^{\otimes d})\simeq
R\Gamma_{et}(\X_{s},\mu_m^{\otimes d}).$$ 
Thus it suffices to show that the cohomology of the right hand side of
$$R\Gamma_{et}(\mathcal{X},\widehat{\Z}'(d))\simeq 
R\mathrm{lim}_{p\nmid m} \,R\Gamma_{et}(\X_{s},\mu_m^{\otimes d})$$
if finite.
By the analog of Proposition \ref{finiteblowup} and induction on the 
number of irreducible components of $\X_{s}$ it suffices to show that the 
cohomology $R\Gamma_{et}(Y,\Z_l(d))$ of each connected component $Y$
of each $\X^{(i)}_{s}$ is finite for all $l\neq p$ and zero for almost all
$l$. Since $Y$ is smooth and proper, this is 
known for the extension $\bar Y$ to the algebraic closure by 
Gabber's theorem \cite{Gabber83}, \cite{Suh},
and this extends to $Y$ by a weight argument because $d>\dim Y$,
hence the Frobenius does not have eigenvalue one on $R\Gamma_{et}(Y,\Z_l(d))$.
\end{proof}

\begin{thm}\label{lemloccompact}
a) Suppose that $(\X_s)^{\red}$ is a simple normal crossing  scheme. Then for any $i\in\Z$, the object $H^{i}_{ar}(\mathcal{X}_{K},\Z)\in \mathcal{LH}(\mathrm{LCA})$ is a discrete abelian group. More precisely, $H^{j}_{ar}(\mathcal{X}_{K},\Z)\in \mathcal{LH}(\mathrm{LCA})$ is an extension of a torsion abelian group by a finitely generated abelian group.

b)  Suppose that $\X/\mathcal{O}_K$ has strictly semi-stable reduction and suppose that $R\Gamma_W(\mathcal{X}_s,\Z^c(0))$ is a perfect complex of abelian groups. Then for any $i\in\Z$, the object  $H^{i}_{ar}(\mathcal{X}_{K},\Z(d))$ is a locally compact abelian group. More precisely, $H^{i}_{ar}(\mathcal{X}_{K},\Z(d))$ is an extension of a finitely generated abelian group by a finitely generated $\Z_p$-module endowed with the $p$-adic topology.
\end{thm}

\begin{proof}
a) We have a long exact sequence in the abelian category  $\mathcal{LH}(\mathrm{LCA})$
$$\cdots\rightarrow 
H^{j}_{ar}(\mathcal{X},\Z)\rightarrow 
H^{j}_{ar}(\mathcal{X}_{K},\Z)\rightarrow 
H^{j+1}_{W}(\mathcal{X}_{s},Ri^{!}\Z)\rightarrow \cdots$$
where $H^{j+1}_{W}(\mathcal{X}_{s},Ri^{!}\Z)
\simeq H^{j}_{W}(\mathcal{X}_{s},Ri^{!}\bq/\Z)$ 
is a discrete torsion abelian group 
(see the proof of Proposition \ref{prop-localduality}) and 
$H^{j}_{ar}(\mathcal{X},\Z)\simeq H^{j}_{W}(\mathcal{X}_s,\Z)$ 
is a discrete finitely generated abelian group by 
Proposition \ref{prop-normalcrossing}. 
Hence $H^{j}_{ar}(\mathcal{X}_{K},\Z)\in \mathcal{LH}(\mathrm{LCA})$ is 
an extension of a torsion abelian group by a finitely generated abelian group. 
It follows that  $H^{j}_{ar}(\mathcal{X}_{K},\Z)\in\mathrm{LCA}$ since 
$\mathrm{LCA}\subset\mathcal{LH}(\mathrm{LCA})$ is stable under extensions 
\cite[Proposition 1.2.29(c)]{Schneiders99}.

b) We have a long exact sequence in $\mathcal{LH}(\mathrm{LCA})$
$$H^{j}_{W}(\mathcal{X}_{s},Ri^{!}\Z(d))\rightarrow 
H^{j}_{ar}(\mathcal{X},\Z(d))\rightarrow 
H^{j}_{ar}(\mathcal{X}_{K},\Z(d))\rightarrow 
H^{j+1}_{W}(\mathcal{X}_{s},Ri^{!}\Z(d))$$
where $H^{j}_{W}(\mathcal{X}_{s},Ri^{!}\Z(d))$ 
is a discrete finitely generated abelian group by assumption.  Moreover, 
$H^{j}_{ar}(\mathcal{X},\Z(d))\in \mathcal{LH}(\mathrm{LCA})$ 
is the group (see Remark \ref{rem-ar-model=et}) 
$$H^{j}_{et}(\mathcal{X},\widehat{\Z}(d))\simeq\prod_{l} 
H^{j}_{et}(\mathcal{X},\Z_l(d)):=
\prod_{l} H^{j}(R\mathrm{lim} 
(R\Gamma_{et}(\mathcal{X},\Z(d))\otimes^{\bL}\Z/l^{\bullet}))$$
which by Lemma \ref{lemfiniteness} is the product of a finite group and the 
finitely generated $\Z_p$-module $H^{j}_{et}(\mathcal{X},\Z_p(d))$ 
endowed with the $p$-adic topology. 
If we can show that the image of the map 
$H^{j}_{W}(\mathcal{X}_{s},Ri^{!}\Z(d))\rightarrow 
H^{j}_{et}(\mathcal{X},\widehat{\Z}(d))$ is finite, then it will follow that 
$H^{j}_{ar}(\mathcal{X}_{K},\Z(d))$ is an extension of a finitely generated abelian group by a profinite abelian group.
Since we have an isomorphism of finitely generated $\Z_p$-modules
$$H^{j}_{W}(\mathcal{X}_{s},Ri^{!}\Z(d))\otimes_{\Z}\Z_p\simeq H^{j}_{et}(\mathcal{X}_{s},Ri^{!}\Z_p(d)),$$
 it is enough to show that the image of the map
$$H^{j}_{et}(\mathcal{X}_{s},Ri^{!}\Z_p(d))\rightarrow H^{j}_{et}(\mathcal{X},\Z_p(d))$$
is finite, or equivalently that the map
$$H^{j}_{et}(\mathcal{X}_{s},Ri^{!}\mathbb{Q}_p(d))\rightarrow 
H^{j}_{et}(\mathcal{X},\mathbb{Q}_p(d))$$
is the zero-map. This follows from Theorem \ref{leminj} by the localization sequence.
\end{proof}

\section{Duality theorems}\label{sectDuality}
The goal of this section is to prove various duality theorems. In particular, 
we prove Theorem \ref{dualityconj} and Corollary \ref{cor-lca} of the introduction. 
Throughout this section, we use the notation and definitions introduced in Section \ref{uncondsect0} and Section \ref{uncondsectd} 
and we assume the following
\begin{Hypothesis}\label{HypDuality}
At least one of the following conditions holds:
\begin{itemize}
\item we have $d\leq 2$;   
\item the scheme $(\X_s)^{\red}$ is a simple normal crossing  scheme and  
$R\Gamma_W(\mathcal{X}_s,\Z^c(0))$ is a perfect complex of abelian groups. 
\end{itemize}
\end{Hypothesis}

In view of Proposition \ref{homology.fin.gen}, $d\leq 2$ implies that $R\Gamma_W(\mathcal{X}_s,\Z^c(0))$ is a perfect complex of abelian groups.

\subsection{Duality with $\Z$-coefficients}\label{sectionDualityZ}

\begin{thm}\label{dualityTHM} Assume Hypothesis \ref{HypDuality}.
Then there is a perfect pairing
$$R\Gamma_{ar}(\mathcal{X}_K,\Z(d))\underline{\otimes}^{\bL} R\Gamma_{ar}(\mathcal{X}_K,\Z) \longrightarrow \Z[-2d]$$
in $\mathbf{D}^b(\mathrm{FLCA})$.
\end{thm}

The rest of Section \ref{sectionDualityZ} is devoted to the proof of Theorem \ref{dualityTHM}. We assume Hypothesis \ref{HypDuality} throughout.

\begin{proof} 
Recall from Proposition \ref{prop-finiteranks} that 
$R\Gamma_{ar}(\mathcal{X}_K,\Z(n))$ belongs to 
$\mathbf{D}^b(\mathrm{FLCA})$ for $n=0,d$, so that the tensor product
$$R\Gamma_{ar}(\mathcal{X}_K,\Z(d))\underline{\otimes}^{\bL} R\Gamma_{ar}(\mathcal{X}_K,\Z)$$
defined in Section \ref{sectionLCA}, makes sense. Moreover, the equivalence
$Ri^{!}\Z^c(0)^{\mathcal{X}}\simeq \Z^c(0)^{\mathcal{X}_s}$
of Remark \ref{rempurityshifted} and the push-forward map $Rf_*\Z^c(0)^{\mathcal{X}_s}\rightarrow \Z^c(0)^s\simeq\Z[0]$
of \cite[Cor.\ 3.2]{Geisser10}
induce trace maps
$$R\Gamma_{W}(\mathcal{X}_s,Ri^{!}\Z(d))\simeq R\Gamma_{W}(\mathcal{X}_s,\Z^c(0)[-2d])\rightarrow R\Gamma_{W}(s,\Z[-2d])\rightarrow \Z[-2d-1]$$
and
$$R\Gamma_{et}(\mathcal{X}_s,Ri^{!}\Z/m(d))\rightarrow R\Gamma_{et}(s,\Z/m[-2d])\rightarrow \Z/m[-2d-1].$$
We start with the following

\begin{prop}\label{Hyp-product}
The canonical product map
$\Z\otimes^{\bL} \Z(d)\rightarrow \Z(d)$ in the derived 
$\infty$-category of \'etale sheaves over $\mathcal{X}_{\mathcal{O}_{K^{un}}}$ and 
$\mathcal{X}_{K^{un}}$ induce perfect pairings
$$R\Gamma_{et}(\mathcal{X}_{s},Ri^!\Z/m)\otimes^{\bL} R\Gamma_{et}(\mathcal{X},\Z/m(d))\rightarrow R\Gamma_{et}(\mathcal{X}_{s},Ri^!\Z/m(d))\rightarrow \Z/m[-2d-1]$$
and
$$R\Gamma_{et}(\mathcal{X}_{s},Ri^!\Z/m(d))\otimes^{\bL} R\Gamma_{et}(\mathcal{X},\Z/m)\rightarrow R\Gamma_{et}(\mathcal{X}_{s},Ri^!\Z/m(d))\rightarrow \Z/m[-2d-1]$$
for any $m$.
\end{prop}

\begin{proof}
Consider the commutative square 
\[ \xymatrix{
R\Gamma_{et}(\mathcal{X},\Z)\otimes^{\bL} R\Gamma_{et}(\mathcal{X},\Z(d))\ar[d]\ar[r]& R\Gamma_{et}(\mathcal{X},\Z(d))\ar[d]\\
R\Gamma_{et}(\mathcal{X}_{K},\Z)\otimes^{\bL} R\Gamma_{et}(\mathcal{X},\Z(d))\ar[r]&R\Gamma_{et}(\mathcal{X}_K,\Z(d)).
}
\]
Taking the fibers of the vertical arrows induces the product map
\begin{equation}\label{pair1}
R\Gamma_{et}(\mathcal{X}_{s},Ri^{!}\Z) \otimes^{\bL} R\Gamma_{et}(\mathcal{X},\Z(d)) \rightarrow R\Gamma_{et}(\mathcal{X}_{s},Ri^{!}\Z(d)),
\end{equation}
and the product map
\begin{equation}\label{pair2}
R\Gamma_{et}(\mathcal{X}_{s},Ri^{!}\Z(d)) \otimes^{\bL} R\Gamma_{et}(\mathcal{X},\Z) \rightarrow R\Gamma_{et}(\mathcal{X}_{s},Ri^{!}\Z(d))
\end{equation}
is obtained similarly.
By \cite[Theorem 7.5]{Geisser10} applied to $\mathcal{F}=\Z/m$, the  pairing
$$R\Gamma_{et}(\mathcal{X}_{s},Ri^!\Z/m)\otimes^{\bL} R\Gamma_{et}(\mathcal{X},\Z/m(d))\rightarrow R\Gamma_{et}(\mathcal{X}_{s},Ri^!\Z/m(d))\rightarrow \Z/m[-2d-1],$$
induced by (\ref{pair1}), is perfect. The pairing induced by (\ref{pair2})
$$R\Gamma_{et}(\mathcal{X}_{s},Ri^!\Z/m(d))\otimes^{\bL} R\Gamma_{et}(\mathcal{X},\Z/m)\rightarrow R\Gamma_{et}(\mathcal{X}_{s},Ri^!\Z/m(d))\rightarrow \Z/m[-2d-1]$$
is perfect as well, since it reduces, by purity and proper base change, to
$$R\Gamma_{et}(\mathcal{X}_{s},\Z^c/m[-2d])\otimes^{\bL} R\Gamma_{et}(\mathcal{X}_s,\Z/m)\rightarrow R\Gamma_{et}(\mathcal{X}_{s},\Z^c/m[-2d])\rightarrow \Z/m[-2d-1]$$
which is perfect by \cite[Theorem 5.1]{Geisser10} applied to $\mathcal{F}=\Z/m$.
\end{proof}

For $n=0$ or $n=d$, consider the product map
$$R\Gamma_{et}(\mathcal{X}_{\bar{s}},R\bar{i}^!\Z(n))\otimes^{\bL} R\Gamma_{et}(\mathcal{X}_{\mathcal{O}_{{K}^{un}}},\Z(d-n))\rightarrow R\Gamma_{et}(\mathcal{X}_{\bar{s}},R\bar{i}^!\Z(d)).$$
This product map is induced by the obvious product maps $\Z\otimes^{\bL} \Z(d)\rightarrow \Z(d)$  in the derived $\infty$-category of \'etale sheaves over $\mathcal{X}_{\mathcal{O}_{K^{un}}}$ and $\mathcal{X}_{K^{un}}$, as in the proof of Proposition \ref{Hyp-product}. Applying $R\Gamma(W_{\kappa(s)},-)$ and composing with the trace map, we obtain
\begin{equation}\label{lili}
R\Gamma_{W}(\mathcal{X}_{s},Ri^!\Z(n))\otimes^{\bL} R\Gamma_{W}(\mathcal{X},\Z(d-n))\rightarrow R\Gamma_{W}(\mathcal{X}_{s},Ri^!\Z(d))\rightarrow \Z[-2d-1].
\end{equation}
This yields the morphisms
$$R\Gamma_{W}(\mathcal{X},\Z(d))\rightarrow R\mathrm{Hom}(R\Gamma_{W}(\mathcal{X}_{s},Ri^!\Z),\Z[-2d-1])$$
which in turn induces
\begin{equation}\label{lala}
\tau^{>a}R\Gamma_{W}(\mathcal{X},\Z(d))\rightarrow R\mathrm{Hom}(R\Gamma_{W}(\mathcal{X}_{s},Ri^!\Z),\Z[-2d-1])
\end{equation}
for $a\ll 0$, since the right hand side is bounded. Composing (\ref{lala}) with the canonical map (see Lemma \ref{lemmapdiscretecont}) 
$$R\mathrm{Hom}(R\Gamma_{W}(\mathcal{X}_{s},Ri^!\Z),\Z[-2d-1])\rightarrow R\underline{\mathrm{Hom}}(R\Gamma_{W}(\mathcal{X}_{s},Ri^!\Z),\Z[-2d-1])$$
 we obtain
\begin{equation}\label{startmap}
\tau^{>a}R\Gamma_{W}(\mathcal{X},\Z(d))\rightarrow R\underline{\mathrm{Hom}}(R\Gamma_{W}(\mathcal{X}_{s},Ri^!\Z),\Z[-2d-1]).
\end{equation}

\begin{prop}\label{prop-localduality}
 The map (\ref{startmap}) factors through an equivalence
\begin{equation}\label{local1}
R\Gamma_{ar}(\mathcal{X},\Z(d))\stackrel{\sim}{\rightarrow} R\underline{\mathrm{Hom}}(R\Gamma_{W}(\mathcal{X}_{s},Ri^!\Z),\Z[-2d-1])
\end{equation}
 in $\mathbf{D}^b(\mathrm{FLCA})$.
\end{prop}

\begin{proof}
One has
\begin{eqnarray*}
R\Gamma_{ar}(\mathcal{X},\Z(d))&:=& \tau^{>a}R\Gamma_{W}(\mathcal{X},\Z(d))\underline{\widehat{\otimes}}\widehat{\Z} \\
&\simeq & (\mathrm{hocolim}\,R\mathrm{Hom}(R\Gamma_{W}(\mathcal{X},\Z/m(d)),\mathbb{Q}/\Z))^D \\
&\simeq&(\mathrm{hocolim}\,R\mathrm{Hom}(R\Gamma_{et}(\mathcal{X},\Z/m(d)),\mathbb{Q}/\Z))^D \\
&\simeq& R\underline{\mathrm{Hom}} (\mathrm{hocolim}\,R\mathrm{Hom}(R\Gamma_{et}(\mathcal{X},\Z/m(d)),\mathbb{Q}/\Z[-2d-1]),\mathbb{R}/\Z[-2d-1])\\
&\stackrel{\sim}{\rightarrow}& R\underline{\mathrm{Hom}} (\mathrm{hocolim}\,R\Gamma_{et}(\mathcal{X}_{s},Ri^!\Z/m),\mathbb{R}/\Z[-2d-1])\\
&\simeq & R\underline{\mathrm{Hom}}(R\Gamma_{W}(\mathcal{X}_{s},Ri^!\mathbb{Q}/\Z),\mathbb{R}/\Z[-2d-1])\\
&\simeq & R\underline{\mathrm{Hom}}(R\Gamma_{W}(\mathcal{X}_{s},Ri^!\Z[1]),\mathbb{R}/\Z[-2d-1])\\
\label{local4}&\simeq& R\underline{\mathrm{Hom}}(R\Gamma_{W}(\mathcal{X}_{s},Ri^!\Z),\Z[-2d-1])
\end{eqnarray*}
where we use Proposition \ref{Hyp-product}, the vanishing
\begin{equation}\label{vanishing}
R\underline{\mathrm{Hom}}(R\Gamma_{W}(\mathcal{X}_{s},Ri^!\Z[1]),\mathbb{R})\simeq 0
\end{equation}
proven in Lemma \ref{lemvanish} below and $Ri^!\mathbb{Q}\simeq 0$. 
\end{proof}
\begin{lem}\label{lemvanish}
We have 
$$R\underline{\mathrm{Hom}}(R\Gamma_{W}(\mathcal{X}_{s},Ri^!\Z[1]),\mathbb{R})\simeq R\underline{\mathrm{Hom}}(\mathbb{R},R\Gamma_{W}(\mathcal{X}_{s},Ri^!\Z[1]))\simeq 0.$$
\end{lem}
\begin{proof}  
As observed above, we have $R\Gamma_{W}(\mathcal{X}_{s},Ri^!\Z[1])\simeq R\Gamma_{W}(\mathcal{X}_{s},Ri^!\mathbb{Q}/\Z)$. Since $R\underline{\mathrm{Hom}}(\mathbb{R},-)$ and $R\underline{\mathrm{Hom}}(-,\mathbb{R})$ are exact functors, and using the $t$-structure on $\mathbf{D}^b(\mathrm{FLCA})$, we may suppose that $R\Gamma_{W}(\mathcal{X}_{s},Ri^!\mathbb{Q}/\Z)$ is cohomologically concentrated in one degree. Hence one is reduced to show that
$$R\underline{\mathrm{Hom}}(A,\mathbb{R})\simeq R\underline{\mathrm{Hom}}(\mathbb{R},A)\simeq 0$$
for any torsion discrete abelian group of finite ranks $A$. This follows from the fact that $\br$ is injective in $\mathrm{LCA}$, see \cite[Proposition 4.15]{Hoffmann-Spitzweck-07}. 
\end{proof}

\begin{cor}\label{lemvanish2}
We have
$$R\underline{\mathrm{Hom}}(R\Gamma_{ar}(\mathcal{X},\Z(d)),\mathbb{R})\simeq R\underline{\mathrm{Hom}}(\mathbb{R},R\Gamma_{ar}(\mathcal{X},\Z(d)))\simeq 0.$$
\end{cor}

\begin{proof}
In the proof of Proposition \ref{prop-localduality}, we have shown that $R\Gamma_{ar}(\mathcal{X},\Z(d))$ is, up to a shift, Pontryagin dual to $R\Gamma_{W}(\mathcal{X}_{s},Ri^!\Z[1])$. Hence the corollary follows from Lemma \ref{lemvanish} since $R\underline{\mathrm{Hom}}(X,Y)\simeq R\underline{\mathrm{Hom}}(Y^D,X^D)$ for any $X,Y\in \mathbf{D}^b(\mathrm{FLCA})$.
\end{proof}

Similarly, we have the

\begin{prop}\label{prop-localduality-dual}
The map
\begin{equation}\label{startmapdual}
R\Gamma_{W}(\mathcal{X}_{s},Ri^{!}\Z(d))\rightarrow R\mathrm{Hom}(R\Gamma_{W}(\mathcal{X},\Z),\Z[-2d-1]),
\end{equation}
induced by (\ref{lili}), factors through an equivalence
\begin{equation}\label{local2}
R\Gamma_{W}(\mathcal{X}_{s},Ri^{!}\Z(d)) \stackrel{\sim}{\longrightarrow}  R\underline{\mathrm{Hom}}(R\Gamma_{ar}(\mathcal{X},\Z) ,\Z[-2d-1]).
\end{equation}
\end{prop}

\begin{proof}
Recall from Remark \ref{rempurityshifted} that we have
$$Ri^!\Z(d)=Ri^!\Z^c(0)[-2d]\simeq \Z^c(0)[-2d].$$
If $\X_s^{\red}$ is a simple normal crossing  scheme, we may therefore identify the map
\begin{equation}\label{local005}
R\Gamma_{W}(\mathcal{X}_{s},Ri^{!}\Z(d))[2d] \stackrel{\sim}{\longrightarrow}  R\underline{\mathrm{Hom}}(R\Gamma_{ar}(\mathcal{X},\Z) ,\Z[-2d-1]) [2d]
\end{equation}
with the composite morphism
\begin{equation*}\label{local5}
R\Gamma_{W}(\mathcal{X}_{s},\Z^c(0))  \stackrel{\sim}{\rightarrow}  R\mathrm{Hom}(R\Gamma_{W}(\mathcal{X}_s,\Z) ,\Z[-1])\stackrel{\sim}{\rightarrow}  R\underline{\mathrm{Hom}}(R\Gamma_{W}(\mathcal{X}_s,\Z) ,\Z[-1])
\end{equation*}
which is an equivalence of perfect complexes of abelian groups by Proposition \ref{prop-normalcrossing}, Theorem \ref{dualityXs}, and Lemma \ref{lemmapdiscretecont}. If $d\leq2$, we may identify (\ref{local005})  with the morphism
\begin{equation*}\label{local50}
R\Gamma_{W}(\mathcal{X}_{s},\Z^c(0))  \stackrel{\sim}{\rightarrow}  R\mathrm{Hom}(R\Gamma_{Wh}(\mathcal{X}_s,\Z) ,\Z[-1])\stackrel{\sim}{\rightarrow}  R\underline{\mathrm{Hom}}(R\Gamma_{Wh}(\mathcal{X}_s,\Z) ,\Z[-1])
\end{equation*}
which is an equivalence of perfect complexes of abelian groups by 
Proposition \ref{propperfectWh} and \cite[Thm.\ 4.2]{Geisser12} 
(using the fact that for a curve, \'etale and eh-cohomology agree).

Note that, if $d\leq 2$, then the following diagram in $\mathbf{D}^b(\mathrm{LCA})$
\[ \xymatrix{
R\Gamma_{W}(\mathcal{X}_{s},Ri^{!}\Z(d)) \ar[rd]^{\simeq} \ar[r]^{(\ref{local2})\hspace{1cm}}\ar[d]^{(\ref{startmapdual})} &  R\underline{\mathrm{Hom}}(R\Gamma_{ar}(\mathcal{X},\Z) ,\Z[-2d-1])\\
R\mathrm{Hom}(R\Gamma_{W}(\mathcal{X},\Z),\Z[-2d-1])  & R\mathrm{Hom}(R\Gamma_{Wh}(\mathcal{X}_s,\Z) ,\Z[-2d-1])\ar[l]\ar[u]^{\simeq}
}
\]
commutes. If $\X_s^{\red}$ is a simple normal crossing  scheme, then the same diagram 
with $R\Gamma_{Wh}(\mathcal{X}_s,\Z)$ replaced by $R\Gamma_{W}(\mathcal{X}_s,\Z)$ commutes as well.

\end{proof}

We now combine Proposition \ref{prop-localduality}  and Proposition \ref{prop-localduality-dual} to prove our result for the generic fiber.

\begin{prop}\label{propfinal}
There is an equivalence $$R\Gamma_{ar}(\mathcal{X}_K,\Z(d))  \stackrel{\sim}{\rightarrow} R\underline{\mathrm{Hom}}(R\Gamma_{ar}(\mathcal{X}_K,\Z) ,\Z[-2d])$$
such that, for any $m$, there is a commutative square
\[ \xymatrix{
R\Gamma_{ar}(\mathcal{X}_K,\Z(d))  \ar[r]^{\sim\hspace{1cm}}\ar[d] &  R\underline{\mathrm{Hom}}(R\Gamma_{ar}(\mathcal{X}_K,\Z) ,\Z[-2d])\ar[d]\\
R\Gamma_{et}(\mathcal{X}_K,\Z/m(d))  \ar[r]^{\sim\hspace{1.5cm}}& R\underline{\mathrm{Hom}}(R\Gamma_{et}(\mathcal{X}_{K},\Z/m),\mathbb{Q}/\Z[-2d]),
}
\]
where the lower horizontal map is induced by duality for the usual \'etale 
cohomology of the variety $\mathcal{X}_K$.
\end{prop}

\begin{proof}
We start with 
the commutative diagram:
\[ \xymatrix{
R\Gamma_{et}(\mathcal{X}_{\bar{s}},R\bar{i}^!\Z)\otimes^{\bL} R\Gamma_{et}(\mathcal{X}_{\mathcal{O}_{{K}^{un}}},\Z(d))\ar[rd]&\\
R\Gamma_{et}(\mathcal{X}_{\bar{s}},R\bar{i}^!\Z)\otimes^{\bL} R\Gamma_{et}(\mathcal{X}_{\bar{s}},R\bar{i}^!\Z(d))\ar[r]\ar[u]\ar[d]& R\Gamma_{et}(\mathcal{X}_{\bar{s}},R\bar{i}^!\Z(d))\\
R\Gamma_{et}(\mathcal{X}_{\mathcal{O}_{{K}^{un}}},\Z)\otimes^{\bL} R\Gamma_{et}(\mathcal{X}_{\bar{s}},R\bar{i}^!\Z(d))\ar[ru]&
}
\]
where the map $$R\Gamma_{et}(\mathcal{X}_{\bar{s}},R\bar{i}^!\Z)\otimes^{\bL} R\Gamma_{et}(\mathcal{X}_{\bar{s}},R\bar{i}^!\Z(d))\rightarrow R\Gamma_{et}(\mathcal{X}_{\bar{s}},R\bar{i}^!\Z(d))$$
is induced by the map $\Z\otimes^{\bL}\Z(d)\rightarrow \Z(d)$ over $\mathcal{X}_{\mathcal{O}_{{K}^{un}}}$ as follows. Consider the morphism
$$\bar{i}_*R\bar{i}^!\Z\rightarrow\Z\rightarrow R\mathbf{Hom}_{\mathcal{X}_{\mathcal{O}_{{K}^{un}}}}(\Z(d),\Z(d))$$
$$\rightarrow R\mathbf{Hom}_{\mathcal{X}_{\mathcal{O}_{{K}^{un}}}}(\bar{i}_*R\bar{i}^!\Z(d),\Z(d))\simeq \bar{i}_*R\mathbf{Hom}_{\mathcal{X}_{\bar{s}}}(R\bar{i}^!\Z(d),R\bar{i}^!\Z(d))$$
and apply $\bar{i}^*$, where $\mathbf{Hom}$ denotes the internal Hom in the category of sheaves on the small \'etale site of the corresponding scheme. Applying $R\Gamma(W_{\kappa(s)},-)$ to the  diagram above, we obtain the following commutative diagram in $\mathbf{D}(\mathrm{Ab})$,
where $tr$ is the trace map:
\[ \xymatrix{
R\Gamma_{W}(\mathcal{X}_{s},Ri^!\Z)\otimes^{\bL} R\Gamma_{W}(\mathcal{X},\Z(d))\ar[rd]& \\
R\Gamma_{W}(\mathcal{X}_{s},Ri^!\Z)\otimes^{\bL} 
R\Gamma_{W}(\mathcal{X}_{s},Ri^!\Z(d))\ar[r]\ar[u]\ar[d]
& R\Gamma_{W}(\mathcal{X}_{s},Ri^!\Z(d)) \ar[r]^{\hspace{.5cm}tr}
& \Z[-2d-1]\\
R\Gamma_{W}(\mathcal{X},\Z)\otimes^{\bL} R\Gamma_{W}(\mathcal{X}_{s},Ri^!\Z(d))\ar[ru]&
}
\]
It gives the following commutative diagram in $\mathbf{D}(\mathrm{Ab})$
\[ \xymatrix{
R\Gamma_{W}(\mathcal{X}_{s},Ri^{!}\Z(d))  \ar[r]\ar[d] &  R\mathrm{Hom}(R\Gamma_{W}(\mathcal{X},\Z) ,\Z[-2d-1])\ar[d]\\
R\Gamma_{W}(\mathcal{X},\Z(d))  \ar[r] \ar[d]& R\mathrm{Hom}(R\Gamma_{W}(\mathcal{X}_{s},Ri^{!}\Z) ,\Z[-2d-1]) \\
\tau^{>a}R\Gamma_{W}(\mathcal{X},\Z(d))\ar[ru]&
}
\]
We obtain the following commutative diagram
\[ \xymatrix{
R\Gamma_{W}(\mathcal{X}_{s},Ri^{!}\Z(d))  \ar[r]^{\eqref{startmapdual}\hspace{1.5cm}}\ar[d] &  R\mathrm{Hom}(R\Gamma_{W}(\mathcal{X},\Z) ,\Z[-2d-1])\ar[d]\\
\tau^{>a}R\Gamma_{W}(\mathcal{X},\Z(d))  \ar[r]^{\eqref{startmap}\hspace{1.5cm}}& R\mathrm{Hom}(R\Gamma_{W}(\mathcal{X}_{s},Ri^{!}\Z) ,\Z[-2d-1]) \ar[d]\\
& R\underline{\mathrm{Hom}}(R\Gamma_{W}(\mathcal{X}_{s},Ri^{!}\Z) ,\Z[-2d-1])
}
\]
in the derived $\infty$-category $\mathbf{D}^b(\mathrm{LCA})$, 
where the lower right map is given by Lemma \ref{lemmapdiscretecont}. 
By construction of the maps (\ref{local1}) and (\ref{local2}), 
we obtain the following commutative diagram
\[ \xymatrix{
&  R\underline{\mathrm{Hom}}(R\Gamma_{ar}(\mathcal{X},\Z) ,\Z[-2d-1])\ar[d]\\
R\Gamma_{W}(\mathcal{X}_{s},Ri^{!}\Z(d))  \ar[r]\ar[d]\ar[ru]^{(\ref{local2})} &  R\mathrm{Hom}(R\Gamma_{W}(\mathcal{X},\Z) ,\Z[-2d-1])\ar[d]\\
\tau^{>a}R\Gamma_{W}(\mathcal{X},\Z(d))  \ar[d]\ar[r]& R\underline{\mathrm{Hom}}(R\Gamma_{W}(\mathcal{X}_{s},Ri^{!}\Z) ,\Z[-2d-1])\\
R\Gamma_{ar}(\mathcal{X},\Z(d))  \ar[ru]_{(\ref{local1})} &
}
\]
hence the upper square in the commutative square
\[ \xymatrix{
R\Gamma_{W}(\mathcal{X}_{s},Ri^{!}\Z(d))  
\ar[r]^{(\ref{local2})\hspace{1.5cm}}_{\sim\hspace{1.5cm}}\ar[d] &  R\underline{\mathrm{Hom}}(R\Gamma_{ar}(\mathcal{X},\Z) ,\Z[-2d-1])\ar[d]\\
R\Gamma_{ar}(\mathcal{X},\Z(d)) 
\ar[r]^{(\ref{local1})\hspace{1.5cm}}_{\sim\hspace{1.5cm}}\ar[d]
& R\underline{\mathrm{Hom}}(R\Gamma_{W}(\mathcal{X}_{s},Ri^{!}\Z) ,\Z[-2d-1])
\ar[d]\\
R\Gamma_{ar}(\mathcal{X}_K,\Z(d))  \ar[r]
& R\underline{\mathrm{Hom}}(R\Gamma_{ar}(\mathcal{X}_K,\Z) ,\Z[-2d]).
}
\]
It follows that the diagram is an equivalence of cofiber sequences in 
$\mathbf{D}^b(\mathrm{LCA})$.
Tensoring the upper commutative square with $\Z/m$ gives a square equivalent to the commutative square
\[ \xymatrix{
R\Gamma_{et}(\mathcal{X}_{s},Ri^{!}\Z/m(d))  \ar[r]\ar[d] &  R\underline{\mathrm{Hom}}(R\Gamma_{et}(\mathcal{X},\Z/m) ,\mathbb{Q}/\Z[-2d-1])\ar[d]\\
R\Gamma_{et}(\mathcal{X},\Z/m(d))  \ar[r]& R\underline{\mathrm{Hom}}(R\Gamma_{et}(\mathcal{X}_{s},Ri^{!}\Z/m) ,\mathbb{Q}/\Z[-2d-1])
}
\]
 where the horizontal maps are induced by the perfect pairings of Proposition \ref{Hyp-product}. This yields the commutative square of Proposition \ref{propfinal}.
\end{proof}

It remains to prove that 
$$R\Gamma_{ar}(\mathcal{X}_K,\Z)\rightarrow R\underline{\mathrm{Hom}}(R\Gamma_{ar}(\mathcal{X}_K,\Z(d)) ,\Z[-2d])$$
is an equivalence.

\begin{lem}\label{lemdual}  
The map
$$R\Gamma_{ar}(\mathcal{X}_K,\Z)\rightarrow R\underline{\mathrm{Hom}}(R\underline{\mathrm{Hom}}(R\Gamma_{ar}(\mathcal{X}_K,\Z) ,\Z),\Z)$$
is an equivalence.
\end {lem}

\begin{proof}
We have 
$$R\Gamma_{ar}(\mathcal{X},\Z)\simeq R\underline{\mathrm{Hom}}(R\underline{\mathrm{Hom}}(R\Gamma_{ar}(\mathcal{X},\Z) ,\Z),\Z)$$
by Lemma \ref{lemmapdiscretecont}, since $R\Gamma_{ar}(\mathcal{X},\Z)$ is a perfect complex of abelian groups by Propositions \ref{prop-normalcrossing}, \ref{propperfectWh}, and \ref{red-map-n=0}. In view of the cofiber sequence
$$R\Gamma_{W}(\mathcal{X}_s,Ri^!\Z)\rightarrow R\Gamma_{ar}(\mathcal{X},\Z)\rightarrow R\Gamma_{ar}(\mathcal{X}_K,\Z)$$
one is reduced to check that the map
$$R\Gamma_{W}(\mathcal{X}_s,Ri^!\Z)\rightarrow R\underline{\mathrm{Hom}}(R\underline{\mathrm{Hom}}(R\Gamma_{W}(\mathcal{X}_s,Ri^!\Z) ,\Z),\Z)$$
is an equivalence. Recall from the proof of Proposition \ref{prop-localduality} that we have 
\begin{eqnarray*}
&&R\underline{\mathrm{Hom}}(R\underline{\mathrm{Hom}}(R\Gamma_{W}(\mathcal{X}_s,Ri^!\Z) ,\Z),\Z)\\
&\simeq& R\underline{\mathrm{Hom}}(R\Gamma_{ar}(\mathcal{X},\Z(d))[2d+1],\Z)\\
&\simeq & R\underline{\mathrm{Hom}}(R\Gamma_{ar}(\mathcal{X},\Z(d))[2d+1],\mathbb{R}/\Z[-1])\\
&\simeq& R\Gamma_{ar}(\mathcal{X},\Z(d))^D[-2d-2]\\
&\simeq & (\mathrm{hocolim}\,R\mathrm{Hom}(R\Gamma_{W}(\mathcal{X},\Z/m(d)),\mathbb{Q}/\Z))^{DD}[-2d-2]\\
&\simeq & \mathrm{hocolim}\,R\mathrm{Hom}(R\Gamma_{W}(\mathcal{X},\Z/m(d)),\mathbb{Q}/\Z[-2d-1])[-1]\\
&\simeq & R\Gamma_{W}(\mathcal{X}_s,Ri^!\mathbb{Q}/\Z)[-1]\\
&\simeq & R\Gamma_{W}(\mathcal{X}_s,Ri^!\Z).
\end{eqnarray*}
where the second equivalence follows from Corollary \ref{lemvanish2}.
\end{proof}
Consider the pairing 
\begin{equation}\label{finalpairing}
R\Gamma_{ar}(\mathcal{X}_K,\Z(d))\underline{\otimes}^{\bL} R\Gamma_{ar}(\mathcal{X}_K,\Z) \longrightarrow \Z[-2d]
\end{equation}
induced by the equivalence of Proposition \ref{propfinal}. The induced map
\begin{equation}\label{finalmap}
R\Gamma_{ar}(\mathcal{X}_K,\Z(d))  \stackrel{\sim}{\rightarrow} R\underline{\mathrm{Hom}}(R\Gamma_{ar}(\mathcal{X}_K,\Z) ,\Z[-2d])
\end{equation}
is (tautologically) the equivalence of Proposition \ref{propfinal}. Moreover, the map
\begin{equation}\label{finalmap2}
R\Gamma_{ar}(\mathcal{X}_K,\Z)\stackrel{\sim}{\rightarrow} R\underline{\mathrm{Hom}}(R\Gamma_{ar}(\mathcal{X}_K,\Z(d)) ,\Z[-2d])
\end{equation}
induced by (\ref{finalpairing}) is an equivalence as well.  Indeed, applying $R\underline{\mathrm{Hom}}(-,\Z[-2d])$ to (\ref{finalmap}) and using Lemma \ref{lemdual}, we obtain the composite equivalence
$$R\Gamma_{ar}(\mathcal{X}_K,\Z)\stackrel{\sim}{\rightarrow} R\underline{\mathrm{Hom}}(R\underline{\mathrm{Hom}}(R\Gamma_{ar}(\mathcal{X}_K,\Z) ,\Z[-2d]),\Z[-2d])$$
$$\stackrel{\sim}{\rightarrow} R\underline{\mathrm{Hom}}(R\Gamma_{ar}(\mathcal{X}_K,\Z(d)),\Z[-2d])$$
which is, up to equivalence, the map (\ref{finalmap2}). 
\end{proof}

\subsection{Pontryagin duality}

Recall that we denote by $\mathrm{FLCA}$ the category of locally compact abelian 
group of finite ranks in the sense of \cite{Hoffmann-Spitzweck-07}. 
It follows from (\ref{tensorFLCA}) and Proposition \ref{prop-finiteranks} that the following definition makes sense. 
\begin{defn} Assume Hypothesis \ref{HypDuality}.
For   $n=0,d$, we define
$$R\Gamma_{ar}(\mathcal{X}_K,\br/\Z(n)):=R\Gamma_{ar}(\mathcal{X}_K,\Z(n))\underline{\otimes}^{\bL}\br/\Z;$$
$$R\Gamma_{ar}(\mathcal{X},\br/\Z(n)):=R\Gamma_{ar}(\mathcal{X},\Z(n))\underline{\otimes}^{\bL}\br/\Z.$$
\end{defn}

\begin{cor}\label{cordualityR/Z}
 Assume Hypothesis \ref{HypDuality}. Then 
one has equivalences
$$R\Gamma_{ar}(\mathcal{X}_{K},\mathbb{R}/\Z)\stackrel{\sim}{\longrightarrow } R\underline{\mathrm{Hom}}(R\Gamma_{ar}(\mathcal{X}_{K},\Z(d)),\mathbb{R}/\Z[-2d])$$
and
$$R\Gamma_{ar}(\mathcal{X}_{K},\Z)\stackrel{\sim}{\longrightarrow } R\underline{\mathrm{Hom}}(R\Gamma_{ar}(\mathcal{X}_{K},\mathbb{R}/\Z(d)),\mathbb{R}/\Z[-2d])$$
in $\mathbf{D}^b(\mathrm{FLCA})$.
\end{cor}
\begin{proof}
By Theorem \ref{dualityTHM} and \cite[Remark 4.3(ii)]{Hoffmann-Spitzweck-07}, we have
\begin{eqnarray*}
R\Gamma_{ar}(\mathcal{X}_K,\Z(d)) &\stackrel{\sim}{\rightarrow} &R\underline{\mathrm{Hom}}(R\Gamma_{ar}(\mathcal{X}_K,\Z) ,\Z[-2d])\\
&\stackrel{\sim}{\rightarrow} &R\underline{\mathrm{Hom}}(R\Gamma_{ar}(\mathcal{X}_K,\Z) , R\underline{\mathrm{Hom}}(\mathbb{R}/\Z,\mathbb{R}/\Z[-2d])\\
&\simeq& R\underline{\mathrm{Hom}}(R\Gamma_{ar}(\mathcal{X}_K,\Z)\underline{\otimes}^{\bL}\mathbb{R}/\Z ,\mathbb{R}/\Z[-2d])\\
&:=&R\underline{\mathrm{Hom}}(R\Gamma_{ar}(\mathcal{X}_K,\mathbb{R}/\Z),\mathbb{R}/\Z[-2d]).
\end{eqnarray*}
Applying the functor $R\underline{\mathrm{Hom}}(-,\mathbb{R}/\Z[-2d])$ and using Pontryagin duality, we obtain the first equivalence of the Corollary.

Similarly, we have
\begin{eqnarray*}
R\Gamma_{ar}(\mathcal{X}_K,\Z) &\stackrel{\sim}{\rightarrow} &R\underline{\mathrm{Hom}}(R\Gamma_{ar}(\mathcal{X}_K,\Z(d)) ,\Z[-2d])\\
&\simeq& R\underline{\mathrm{Hom}}(R\Gamma_{ar}(\mathcal{X}_K,\Z(d))\underline{\otimes}^{\bL}\mathbb{R}/\Z ,\mathbb{R}/\Z[-2d])\\
&:=&R\underline{\mathrm{Hom}}(R\Gamma_{ar}(\mathcal{X}_K,\mathbb{R}/\Z(d)),\mathbb{R}/\Z[-2d]).
\end{eqnarray*}

\end{proof}


\begin{cor}
Suppose that $\X/\mathcal{O}_K$ has strictly semi-stable reduction and suppose that $R\Gamma_W(\mathcal{X}_s,\Z^c(0))$ is a perfect complex of abelian groups.
Then for any $i\in\Z$, we have an isomorphism of locally compact groups
$$H^{i}_{ar}(\mathcal{X}_{K},\mathbb{R}/\Z)\stackrel{\sim}{\longrightarrow } H^{2d-i}_{ar}(\mathcal{X}_{K},\Z(d))^D$$
and an isomorphism of discrete groups
$$H^{i}_{ar}(\mathcal{X}_{K},\Z)\stackrel{\sim}{\longrightarrow } H^{2d-i}_{ar}(\mathcal{X}_{K},\mathbb{R}/\Z(d))^D.$$
\end{cor}
\begin{proof}
In view of  Theorem \ref{lemloccompact} and Lemma \ref{remcohdual}, the equivalence in $\mathbf{D}^b(\mathrm{FLCA})$
$$R\Gamma_{ar}(\mathcal{X}_{K},\mathbb{R}/\Z)\stackrel{\sim}{\longrightarrow } R\Gamma(\mathcal{X}_{K},\Z(d))^D[-2d]$$
induces isomorphisms
$$H^{i}_{ar}(\mathcal{X}_{K},\mathbb{R}/\Z)\stackrel{\sim}{\longrightarrow } H^{i}(R\Gamma(\mathcal{X}_{K},\Z(d))^D[-2d])\simeq H^{2d-i}_{ar}(\mathcal{X}_{K},\Z(d))^D$$
of locally compact abelian groups. Similarly, the equivalence in $\mathbf{D}^b(\mathrm{FLCA})$
$$R\Gamma_{ar}(\mathcal{X}_{K},\mathbb{R}/\Z(d))\stackrel{\sim}{\longrightarrow } R\Gamma(\mathcal{X}_{K},\Z)^D[-2d]$$
induces isomorphisms
$$H^{i}_{ar}(\mathcal{X}_{K},\mathbb{R}/\Z(d))\stackrel{\sim}{\longrightarrow } H^{i}(R\Gamma(\mathcal{X}_{K},\Z)^D[-2d])\simeq H^{2d-i}_{ar}(\mathcal{X}_{K},\Z)^D$$
of compact abelian groups.

\end{proof}

\section{The conjectural picture}
We conjecture the existence of a cohomology theory on the category of separated 
schemes of finite type over $\mathrm{Spec}(\mathcal{O}_K)$, 
with values in $\mathbf{D}^b(\mathrm{FLCA})$, which we denote by
$$R\Gamma_{ar}(-,A(n))$$
for any $A\in\mathrm{FLCA}$ and any $n\in\Z$. Furthermore we conjecture that
the conclusion of Theorem \ref{dualityconj} holds in full generality: 
For any smooth proper $\X_K$ over $K$ of pure dimension $d-1$, any
Tate twist $n\in\Z$, and any $A\in\mathrm{FLCA}$, there is an  equivalence
$$R\Gamma_{ar}(\X_K,A^D(n))\stackrel{\sim}{\longrightarrow} 
R\underline{\mathrm{Hom}}(R\Gamma_{ar}(\X_K,A(d-n)), \mathbb{R}/\Z[-2d])$$
induced by a trace map 
$H^{2d}_{ar}(\X_K,\mathbb{R}/\Z(d))\rightarrow \mathbb{R}/\Z$, 
where $A^D$ denotes the Pontryagin dual of $A$. 
However, we do not expect the analog of Corollary \ref{cor-lca} to 
be true in general, since the
groups $H^{i}_{ar}(\X_K,\Z(n))$ cannot be expected to be locally compact 
for arbitrary Tate twist $n$ 
(as one can see from \cite{Geisser-Morin20} for $n=1$). 
Instead, they could be seen as condensed abelian groups 
(or in the language used in this paper, as objects of the heart 
$\mathcal{LH}(\mathrm{LCA})$ of the left $t$-structure on 
$\mathbf{D}^b(\mathrm{LCA})$). 
In contrast, we do expect isomorphisms of \emph{compact} abelian groups 
$$H^{i}_{ar}(\X,A(n))\simeq H^{2d+1-i}_{ar}(\X_s,Ri^!A^D(d-n))^D$$
for any $i,n\in\Z$ and any \emph{compact} $A\in\mathrm{FLCA}$. 
Concerning the relationship between  $R\Gamma_{ar}(-,A(n))$
and  known cohomology theories, we expect the following, for
$\X$ a regular proper flat scheme over $\mathrm{Spec}(\mathcal{O}_K)$.

\begin{itemize}

\item For any positive integer $m$,  we have
$$R\Gamma_{ar}(\X_K,\Z/m(n)) \simeq R\Gamma_{et}(\X_K,\Z/m(n)).$$
In particular, for any prime $l$, one has equivalences
$$ R\Gamma_{ar}(\X_K,\Z_l(n))\simeq R\Gamma_{ar}(\X_K,\Z(n))\underline{\widehat{\otimes}}\Z_l\simeq R\Gamma_{et}(\X_K,\Z_l(n))$$
where $(-)\underline{\widehat{\otimes}}\Z_l:=R\mathrm{lim}(-\underline{\otimes}^{\bL}\Z/l^{\bullet})$ is the $l$-adic completion functor. 

\item  The canonical map
$$R\Gamma_{ar}(-,\Z(n))\underline{\otimes}^{\bL} A\stackrel{\sim}{\rightarrow} R\Gamma_{ar}(-,A(n))$$
is an equivalence for $(-)=\X,\X_s$ and any ring object $A$, and for $(-)=\X_K$ if $A$ has no topological $p$-torsion\footnote{The locally compact group $A$ has a unique filtration by closed subgroups with graded pieces $A_{\mathbb{S}^1}$, $A_{\mathbb{A}}$ and $A_{\Z}$ of type $\mathbb{S}^1$, $\mathbb{A}$ and $\Z$ respectively. Then $A_{\mathbb{A}}$ is the direct sum of a finite dimensional $\br$-vector space and topological torsion group $A_{\mathrm{toptor}}$, which in turn has a topological $p$-torsion component $A_p$ (see \cite[Section 2]{Hoffmann-Spitzweck-07}). We say that  $A$ has no topological $p$-torsion if $A_p=0$.}. For example the map
$$R\Gamma_{ar}(\X_K,\Z(n))\underline{\otimes}^{\bL} A\stackrel{\sim}{\rightarrow} R\Gamma_{ar}(\X_K,A(n))$$
is an equivalence for $A=\br$ and $A=\bq_l$ if $l\neq p$. However, one has an equivalence
$$
R\Gamma_{ar}(\X_K,\Z(n))\underline{\otimes}^{\bL} \bq_p\simeq R\Gamma_{syn}(\X_K,n)
$$
where the right hand side is the Nekovar-Niziol syntomic cohomology \cite{Nekovar-Niziol}. The induced map
$$
R\Gamma_{syn}(\X_K,n)\simeq R\Gamma_{ar}(\X_K,\Z(n))\underline{\otimes}^{\bL}\mathbb{Q}_p \rightarrow R\Gamma_{ar}(\X_K,\Z(n))\underline{\widehat{\otimes}}\mathbb{Q}_p\simeq R\Gamma_{et}(\X_K,\mathbb{Q}_p(n))
$$
is expected to be an equivalence if and only if $n\geq d$.

\item We have equivalences
$$R\Gamma_{ar}(\X,\Z_p(n))\simeq R\Gamma_{ar}(\X,\Z(n))\underline{\widehat{\otimes}}\Z_p\simeq R\Gamma_{et}(\X,\Z(n))\underline{\widehat{\otimes}}\Z_p$$
where $R\Gamma_{et}(\X,\Z(n))$ denotes \'etale motivic cohomology, i.e. \'etale hypercohomology of Bloch's motivic complex. Note that $R\Gamma_{et}(\X,\Z(n))\widehat{\underline{\otimes}}\mathbb{Q}_p$ is equivalent to the  syntomic cohomology of Fontaine-Messing \cite{Fontaine-Messing}, at least if $\X/\mathcal{O}_K$ is smooth and $n<p-1$ (see \cite[Proposition 7.21, Remark 7.23]{Flach-Morin-17}).  For general regular proper flat $\X$ and arbitrary $n$, a conjectural syntomic description of $R\Gamma_{et}(\X,\Z(n))\widehat{\underline{\otimes}}\mathbb{Q}_p$ and of the fiber sequence $(\ref{IntroTri})\widehat{\underline{\otimes}}\mathbb{Q}_p$ is given by \cite[Corollary 7.17]{Flach-Morin-17}.
 
\item We have $$R\Gamma_{ar}(\X_s,\Z(n))\simeq R\Gamma_{Wh}(\X_s,\Z(n)) $$
where $R\Gamma_{Wh}(\X_s,\Z(n))$ is motivic $Wh$-cohomology \cite{Geisser06}.
Moreover,  the cofiber 
$$C_{ar}(\X,n):=\mathrm{Cofib} 
\left(  R\Gamma_{ar}(\X,\Z(n))\rightarrow R\Gamma_{ar}(\X_s,\Z(n))\right)$$ 
is a perfect complex of $\Z_p$-modules such that
$$C_{ar}(\X,n)\otimes\bq\simeq R\Gamma(\X_K,\Omega_{\X_K/K}^{<n})$$
where the right hand side denotes de Rham cohomology modulo the $n$-step of the Hodge filtration.

\item 
We have  
$$\mathrm{dim}_{\bq_l} H_{ar}^{i}(\X_K,\bq_l(n))=\mathrm{dim}_{\br} H_{ar}^{i}(\X_K,\br(n))$$
for any $i,n\in\Z$ and any  prime $l\neq p$. In particular, the left hand side is independent on $l\neq p$.

\end{itemize}

\end{document}